\documentclass[graybox]{svmult}

\usepackage{type1cm}        
\usepackage{makeidx}         
\usepackage{graphicx}        
\usepackage{multicol}        
\usepackage[bottom]{footmisc}

\usepackage{amssymb,epsfig,amsmath,accents,setspace}

\spdefaulttheorem{assumption}{Assumption}{\bfseries }{}

\usepackage[sort&compress]{natbib}
\bibpunct{(}{)}{,}{a}{,}{,}

\newcommand{\bbR}{{\mathbb R}}

\newcommand{\bbX}{{\mathbb X}}

\newcommand{\probP}{\mathbb{P}}

\newcommand{\expect}{\mathbb{E}}
\newcommand{\vari}{\mathrm{var}}

\newcommand{\tran}{^\top}

\newcommand{\trace}{\mathrm{tr}}

\makeindex             

\begin{document}

\title*{Who Are I: Time Inconsistency and Intrapersonal Conflict and Reconciliation
}

\author{Xue Dong He and Xun Yu Zhou}
\institute{Xue Dong He \at The Chinese University of Hong Kong, \email{xdhe@se.cuhk.edu.hk}
\and Xun Yu Zhou \at Columbia University, \email{xz2574@columbia.edu}}
%
%
\maketitle


\abstract{Time inconsistency is prevalent in dynamic choice problems: a plan of actions to be taken in the future that is optimal for an agent today may not be optimal for the same agent in the future. If the agent is aware of this intrapersonal conflict  but unable to commit herself in the future to following the optimal plan today, the rational strategy for her today is to reconcile with her future selves, namely to correctly anticipate her actions in the future and then act today accordingly. Such a strategy is named intra-personal equilibrium and has been studied since  as early as in the 1950s. A rigorous treatment in continuous-time settings, however, had  not been available until a decade ago. Since then, the study on intra-personal equilibrium for time-inconsistent problems in continuous time has grown rapidly. In this chapter, we review the classical results and some recent development in this literature.}

\section{Introduction}
\label{sec:Intro}
When making dynamic decisions, the decision criteria of an agent at different times may not align with each other, leading to time-inconsistent behavior: an action that is optimal under the decision criterion today may no longer be optimal under the decision criterion at certain future time. A variety of preference models can lead to time inconsistent behaviors, such as those involving  present-bias, mean-variance criterion, and probability weighting. 

In his seminal paper, \cite{Strotz1955:MyopiaInconsistency} describes three types of agents when facing time inconsistency. Type 1, a ``spendthrift" (or a naivet\'e as in the more recent literature), does not recognize the time-inconsistency and at any given time seeks an optimal solution from the vantage point of that moment only. As a result, his strategies are always myopic and change all the times. The next two types are aware of the time inconsistency but act differently. Type 2 is a ``precommitter" who solves the optimization problem only once at time 0 and then commits to the resulting strategy throughout, even though she knows that the original solution may no longer be optimal at later times. Type 3 is a ``thrift" (or a sophisticated agent) who is unable to precommit and realizes that her future selves may disobey whatever plans she makes now. Her resolution is to compromise and choose {\it consistent planning} in the sense that she optimizes taking the future disobedience as a {\it constraint}. In this resolution, the agent's selves at different times are considered to be the players of a game, and a consistent plan chosen by the agent becomes an equilibrium of the game from which no selves are willing to deviate. Such a plan or strategy
is referred to as an {\em intra-personal equilibrium}.

To illustrate the above three types of behavior under time inconsistency, consider an agent who has a planning horizon with a finite end date $T$ and makes decisions at {\em discrete} times $t\in \{0,1,\dots, T-1\}$. The agent's decision drives a Markov state process and the agent's decision criterion at time $t$ is to maximize an objective function $J(t,x;\mathbf{u})$, where $x$ stands for the Markovian state at that time and $\mathbf{u}$ represents the agent's strategy. The agent considers Markovian strategies, so $\mathbf{u}$ is a function of time $s\in\{0,1,\dots, T-1\}$ and the Markovian state at that time. If the agent, at certain time $t$ with state $x$, is a ``pre-committer", she is committed to implementing throughout the remaining horizon the strategy $\mathbf{u}^{\mathrm{pc}}_{(t,x)}=\{\mathbf{u}^{\mathrm{pc}}_{(t,x)}(s,\cdot)|s=t,t+1,\dots, T-1\}$ that maximizes $J(t,x;\mathbf{u})$, and this strategy is referred to as the {\em pre-committed strategy} of the agent at time $t$ with state $x$. If the agent is a ``spendthrift", at {\em every} time $t$ with state $x$, she is able to implement the pre-committed strategy at that moment only and will change at the next moment; so the strategy that is actually implemented by the agent throughout the horizon is $\mathbf{u}^{\mathrm{n}}=\{\mathbf{u}^{\mathrm{pc}}_{(s,X^{\mathbf{u}^{\mathrm{n}}}(s))}(s,X^{\mathbf{u}^{\mathrm{n}}}(s))|s=0,1,\dots, T-1\}$, where $X^{\mathbf{u}^{\mathrm{n}}}$ denotes the state process under $\mathbf{u}^{\mathrm{n}}$. This strategy is referred to as the {\em na\"ive strategy}. If the agent is a ``thrift", she chooses an intra-personal equilibrium strategy $\hat{\mathbf{u}}$: At any time $t\in \{0,1,\dots, T-1\}$ with any state $x$ at that time, $\hat{\mathbf{u}}(t,x)$ is the optimal action of the agent given that her future selves follow $\hat{\mathbf{u}}$; i.e.,
\begin{align}\label{eq:EquilibriumDiscreteTime}
 \hat{\mathbf{u}}(t,x) \in \mathrm{arg}\max_{u}J(t,x;  {\mathbf{u}}_{t,u}),
\end{align}
where ${\mathbf{u}}_{t,u}(t,x):=u$  and ${\mathbf{u}}_{t,u}(s,\cdot):=\hat{\mathbf{u}}(s,\cdot)$ for $s=t+1,\dots, T-1$.

All the three types of behavior are important from an economic perspective. First, field and experimental studies reveal the popularity of commitment devices to help individuals to fulfill plans that would otherwise be difficult to implement due to lack of self control; see for instance \citet{BryanEtal2010:CommitmentDevices}. The demand for commitment devices implies that some individuals seek for pre-committed strategies in the presence of time inconsistency. Second, empirically observed decision-making behavior implies that some individuals are naivet\'es. For example, \citet{Barberis2012:Casino} shows that a na\"ive agent would take on a series independent, unfavorable bets and take a gain-exit strategy, and this gambling behavior is commonly observed in casinos. Finally, when an agent foresees the time-inconsistency and a commitment device is either unhelpful or unavailable, the intra-personal equilibrium strategy becomes a rational choice of the agent.

It is important to note that it is hard or perhaps not meaningful to determine which type is superior than the others, simply because there is no uniform criteria to evaluate  and compare them.
So a na\"ive strategy, despite its name,  is not necessarily inferior to an intra-personal equilibrium in terms of an agent's long-run utility.
Indeed, \citet{ODonogheRabin:1999DoingItNowOrLater} show that in an optimal stopping problem with an immediate reward and present-biased preferences, a sophisticate agent has a larger tendency to preproperate than a naivet\'e and thus leads to a lower long-run utility. In this sense, studying the different behaviors under time inconsistency sometimes falls into the realm of being ``descriptive" as in behavioral science, rather than being ``normative" as in classical decision-making theory.

In this survey article, we focus on reviewing the studies on intra-personal equilibrium of a sophisticated agent in {\em continuous time}.\footnote{Hence the title of this article. Note that there is no grammatical error in the phrase ``who {\it are} I": the word ``I" here is {\it plural}; it refers to many different selves over time, which is the premise of this article.}
Intra-personal equilibrium for time-inconsistent problems in discrete time, which is defined through the equilibrium condition \eqref{eq:EquilibriumDiscreteTime},  has been extensively studied in the literature and generated various economic implications. The extension to the continuous-time setting, however, is nontrivial because in this setting, taking a different action from a given strategy at only one time instant  does not change the state process and thus has no impact on the objective function value. As a result, it becomes meaningless to examine whether the agent is willing to deviate from a given strategy at a particular moment  by just comparing the objective function values before and after the deviation. To address this issue and to  formalize the idea of \cite{Strotz1955:MyopiaInconsistency}, \cite{EkelandPirvu2008:InvestConsptnWithoutCommit}, \cite{EkelandLazrak2006:BeingSeriousAboutNonCommitment}, and \cite{BjordTMurgoci:08u} assume that the agent's self at each time can implement her strategy in an infinitesimally small, but positive, time period; consequently, her action has an impact on the state process and thus on the objective function. In Section \ref{sec:ExHJB} below, we follow the framework of \cite{BjorkMurgoci2014:TheoryMarkovianTimeInconsistent} to define intra-personal equilibria, show a sufficient and necessary condition for an equilibrium, and present the so-called {\em extended HJB equation} that characterizes the intra-personal equilibrium strategy and the  value under this strategy. In Section \ref{sec:Discussions}, we further discuss various issues related to intra-personal equilibria.

A close-loop strategy for a control system is a mapping from the historical path of the system state and control to the space of controls. So at each time, the control taken by the agent is obtained by plugging the historical path into this mapping. For example, a Markovian strategy for a Markovian control system is a closed-loop strategy. An open-loop strategy is a collection of controls across time (and across scenarios in case of stochastic control), and at each time the control in this collection is taken, regardless of the historical path of the system state and control. For a classical, time-consistent controlled Markov decision problem, the optimal close-loop strategy and the optimal open-loop strategy yield the same state-control path. For time-inconsistent problems, however, closed-loop and open-loop intra-personal equilibria can be vastly different. In Section \ref{sec:OpenLoop}, we review the study of open-loop intra-personal equilibrium and discuss its connection with closed-loop intra-personal equilibrium.

Optimal stopping problems can be viewed as a special case of control problems, so intra-personal equilibria can be defined similarly for time-inconsistent stopping problems. These problems, however, have very special structures, and by exploiting these structures new notions of intra-personal equilibria have been proposed in the literature. We discuss these in Section \ref{sec:Stopping}.

If we discretize a continuous horizon of time and assume that the agent has full self control in each subperiod under the discretization, we can define and derive intra-personal equilibria as in the discrete-time setting. The limits of the intra-personal equilibria as discretization becomes infinitely finer  are used by some authors to define  intra-personal equilibria for continuous-time problems. In Section \ref{sec:Discretization}, we review this thread of research.

Time-inconsistency arises in various economic problems, and for many of them, intra-personal equilibria have been studied and their implications  discussed in the literature. In Section \ref{sec:Applications}, we review this literature.

Finally, in Section \ref{sec:DynamicConsistent}, we review the studies on dynamic consistency preferences. In these studies, starting from a preference model for an agent at certain initial time, the authors attempt to find certain preference models for the agent's future selves such that the pre-committed strategy for the agent at the initial time is also optimal for the agent at any future time and thus can be implemented consistently over time.

\section{Extended HJB Equation}
\label{sec:ExHJB}

\citet{Strotz1955:MyopiaInconsistency} is the first  to study the behavior of a sophisticated agent in the presence of time-inconsistency in a continuous-time model. Without formally defining the notion of intra-personal equilibrium, the author derives a consistent plan of the sophisticated agent. \citet{Barro1999:Ramsey} and \citet{LuttmerMariotti2003:SubjectiveDiscounting} also investigate, for certain continuous-time models, consistent plans of sophisticated agents, again without their formal definitions. In a series of papers, \citet{EkelandLazrak2006:BeingSeriousAboutNonCommitment}, \citet{EkelandILazrakA:07ti}, and \citet{EkelandLazrak2010:GoldenRule} study the classical Ramsey model with a nonexponential discount function and propose for the first time a formal notion of intra-personal equilibrium for deterministic control problems in continuous time. Such a notion is proposed in a stochastic context by \citet{BjordTMurgoci:08u}, which is later split into two papers, \cite{BjorkMurgoci2014:TheoryMarkovianTimeInconsistent} and \citet{Bjork2017:TimeInconsistent}, discussing the discrete-time and continuous-time settings, respectively. In this section, we follow the framework of \citet{Bjork2017:TimeInconsistent} to define an intra-personal equilibrium strategy and present a sufficient and necessary condition for such a strategy.

\subsection{Notations}

We first introduce some notations. By convention, $x\in\mathbb{R}^n$ is always a column vector. When a vector $x$ is a row vector, we write it as $x\in \mathbb{R}^{1\times n}$. Denote by $A\tran $ the transpose of a matrix $A$, and by $\trace(A)$ the trace of a square matrix $A$. For a differentiable function $\xi$ that maps $x\in \mathbb{R}^m$ to $\xi(x)\in \mathbb{R}^n$, its derivative, denoted as $\xi_x(x)$, is an $n\times m$ matrix with the entry in the $i$-th row and $j$-th column denoting the derivative of the $i$-th component of $\xi$ with respect to the $j$-th component of $x$. In particular, for a mapping $\xi$ from $\mathbb{R}^m$ to $\mathbb{R}$, $\xi_x(x)$ is an $m$-dimensional row vector, and we further denote by $\xi_{xx}$ the Hessian matrix.

Consider $\xi$ that maps $(z,x)\in \mathbb{Z}\times \bbX$ to $\xi(z,x)\in\bbR^l$, where $\mathbb{Z}$ is a certain set  and $\bbX$, which represents the state space throughout, is either $\bbR^n$ or $(0,+\infty)$. $\xi$ is {\em locally Lipschitz in $x\in \bbX$, uniformly in $z\in\mathbb{Z}$} if there exists a sequence of compact sets $\{\bbX_k\}_{k\ge 1}$ with $\cup_{k\ge 1}\bbX_k=\bbX$ and a sequence of positive numbers $\{L_k\}_{k\ge 1}$ such that for any $k\ge 1$, $\|\xi(z,x)-\xi(z,x') \|\leq L_k \|x-x'\|,\forall z\in \mathbb{Z}, x,x'\in \bbX_k$.
 $\xi$ is {\em global Lipschitz in $x\in \bbX$, uniformly in $z\in\mathbb{Z}$} if there exists  constant  $L>0$ such that  $\|\xi(z,x)-\xi(z,x') \|\leq L  \|x-x'\|,\forall z\in \mathbb{Z}, x,x'\in \bbX$. In the case $\bbX=\bbR^n$, $\xi$ is of {\em linear growth in $x\in \bbX$, uniformly in $z\in\mathbb{Z}$} if there exists $L>0$ such that $\|\xi(z,x)\|\le L(1+\|x\|),\forall z\in \mathbb{Z}, x\in \bbX$. In the case $\bbX=(0,+\infty)$, $\xi$ {\em has a bounded norm in $x\in\bbX$, uniformly in $z\in\mathbb{Z}$}, if there exists $L>0$ such that $\|\xi(z,x)\|\le Lx,\forall z\in \mathbb{Z}, x\in \bbX$.
 $\xi$ is of {\em polynomial growth in $x\in\mathbb{X}$, uniformly in $z\in\mathbb{Z}$} if there exists $L>0$ and integer $\gamma\ge 1$ such that $\|\xi(z,x)\|\le L\left(1+\varphi_{2\gamma}(x)\right), \forall z\in \mathbb{Z}, x\in \bbX$, where $\varphi_{2\gamma}(x)=\|x\|^{2\gamma}$ when $\bbX=\bbR^n$ and $\varphi_{2\gamma}(x) = x^{2\gamma}+x^{-2\gamma}$ when $\bbX=(0,+\infty)$.


Fix integers $r\geq 0$, $q\geq 2 r$, and real numbers $a<b$. Consider ${\xi}$ that maps $(t,x)\in [a,b]\times \bbX$ to ${\xi}(t,x) \in \bbR^l $. We say $\xi\in \mathfrak{C}^{r,q}([a, b]\times \bbX )$ if for any derivative index $\alpha$ with $|\alpha|\le q-2j$ and $j=0,\dots, r$, the partial derivative $\frac{\partial^{j+\alpha} \xi(t,x) }{\partial t^j \partial x^\alpha}:=\frac{\partial^{j+\alpha_1+\dots+\alpha_n} \xi(t,x) }{\partial t^j \partial x_1^{\alpha_1}\dots \partial x_n^{\alpha_n}} $ exists for any $(t, x)\in ( a,b) \times \bbX$ and can be extended to and continuous on $[a, b]\times \bbX$.
We say ${\xi}\in \mathfrak{\bar C}^{r,q}( [a, b] \times \bbX)$ if ${\xi}\in \mathfrak{C}^{r,q}( [a, b] \times \bbX)$ and $\frac{\partial^{j+\alpha} \xi(t,x) }{\partial t^j \partial x^\alpha}$ is of polynomial growth in $x\in \bbX$, uniformly in $t\in [a, b]$, for any derivative index $\alpha$ with $|\alpha|\le q-2j$ and $j=0,...,r$. 

\subsection{Time-Inconsistent Stochastic Control Problems}\label{subse:framework}
Let be given a probability space $(\Omega,{\cal F},\probP)$ with a standard $d$-dimensional Brownian motion $W(t):=\big(W_1(t),...,W_d(t)\big)\tran$, $t\ge 0$, on the space, along with the filtration  $({\cal F}_t)_{t\ge 0}$ generated by the Brownian motion and augmented by the $\probP$-null sets.
Consider an agent who makes dynamic decisions in a given period $[0,T]$, and for any $(t,x)\in[0,T)\times \bbX$, the agent faces the following stochastic control problem:
\begin{align}\label{eq:ControlProblem}
\left\{
\begin{array}{c l}
\underset{\mathbf{u}}{\max}  &J(t,x;\mathbf{u} )  	 \\
\text{subject to} &dX^{\mathbf{u}}(s)=\mu(s, X^{\mathbf{u}}(s), \mathbf{u}(s, X^{\mathbf{u}}(s)) )ds\\
 & \;\;\;\;  +\sigma(s,X^{\mathbf{u}}(s), \mathbf{u}(s,X^{\mathbf{u}}(s)) )dW(s),\; s\in[t,T]\\
 & X^{\mathbf{u}}(t)=x.
 \end{array}\right.
 \end{align}
The agent's dynamic decisions are represented by a Markov strategy $\mathbf{u}$, which maps $(s,y)\in [0,T)\times \bbX$ to $\mathbf{u}(s,y)\in \mathbb{U}\subset\bbR^m$. 
The controlled diffusion process $X^{\mathbf{u}}$ under $\mathbf{u}$ takes values in $\bbX$, which as aforementioned is set to be either $(0,+\infty)$ or $\bbR^n$.
$\mu$ and $\sigma$ are measurable mappings from $[0,T]\times \bbX\times \mathbb{U}$ to $\bbR^n$ and to $\bbR^{n\times d}$, respectively, where $n$ stands for the dimension of $\bbX$.

The agent's goal at $(t,x)\in [0,T]\times \bbX$ is to maximize the following objective function:
\begin{align}\label{eq:ObjFun}
J(t,x;\mathbf{u} ) &=\expect_{t,x}\left[\int_t^T C\big(t,x, s, X^{\mathbf{u}}(s), {\mathbf{u}}(s, X^{\mathbf{u}}(s) )    \big)ds+F\big(t, x, X^{{\mathbf{u}}}(T) \big)\right]\notag \\
 &\quad +G\big(t, x, \expect_{t,x} [  X^{{\mathbf{u}}}(T) ]  \big),
\end{align}
where $C$ is a measurable mapping from $[0,T)\times \bbX\times [0,T]\times \bbX\times \mathbb{U}$ to $\mathbb{R}$, and $F$ and $G$ are measurable mappings from $[0,T)\times \bbX\times \bbX$ to $\bbR$. Here and hereafter, $\expect_{t,x}[Z]$ denotes the expectation of $Z$ conditional on $X^{\mathbf{u}}(t)=x$. If $C$, $F$, and $G$ are independent of $(t,x)$ and $G\big(t, x, \expect_{t,x} [  X^{{\mathbf{u}}}(T) ]\big)$ is linear in $\expect_{t,x} [  X^{{\mathbf{u}}}(T) ]$, then $J(t,x;\mathbf{u} )$ becomes a standard objective function in classical stochastic control  where time consistency holds. Thus, with objective function
\eqref{eq:ObjFun}, time inconsistency arises from the dependence of $C$, $F$, and $G$ on $(t,x)$ as well as from  the nonlinearity of $G\big(t, x, \expect_{t,x} [  X^{{\mathbf{u}}}(T) ]\big)$ in $\expect_{t,x} [  X^{{\mathbf{u}}}(T) ]$.

For any feedback strategy $\mathbf{u}$, denote
\begin{align*}
&\mu^{\mathbf{u}}(t, x):=\mu(t, x, \mathbf{u}(t,x)), \; \sigma^{\mathbf{u}}(t, x):=\sigma(t, x, \mathbf{u}(t,x)), \; \nonumber\\
&\Upsilon^{\mathbf{u}}(t, x):=\sigma(t, x, \mathbf{u}(t,x)) \sigma(t, x, \mathbf{u}(t,x))\tran , \; C^{\tau, y, \mathbf{u}}(t,x ):=C\big(\tau, y,t, x, \mathbf{u}(t, x ) \big).
\end{align*}
With a slight abuse of notation, $u \in \mathbb{U}$ also denotes the feedback strategy $\mathbf{u}$ such that $\mathbf{u}(t,x)=u,\forall (t,x)\in[0,T]\times \bbX$; so $\mathbb{U}$ also stands for the set of all {\em constant} strategies when no ambiguity arises.

We need to impose conditions on a strategy $\mathbf{u}$ to ensure the existence and uniqueness of the SDE in \eqref{eq:ControlProblem} and the well-posedness of the objective function $J(t,x;\mathbf{u} )$. This consideration leads to the following definition of feasibility:
\begin{definition}\label{de:Feasibility}
	A feedback strategy $\mathbf{u}$ is {\em feasible} if the following hold:
	\begin{enumerate}
\item[(i)] $\mu^{\mathbf{u}}$, $\sigma^{\mathbf{u}}$ are locally Lipschitz in $x\in\bbX$, uniformly in $t\in[0,T]$.
\item[(ii)] $\mu^{\mathbf{u}}$ and $\sigma^{\mathbf{u}}$ are of linear growth in $x\in\bbX$, uniformly in $t\in[0,T]$, when $\bbX=\bbR^n$ and have bounded norm in $x\in\bbX$, uniformly in $t\in[0,T]$, when $\bbX=(0,+\infty)$.
\item[(iii)] For each fixed $(\tau,y) \in [0, T)\times \bbX$, $C^{\tau,y,\mathbf{u}}(t,x)$ and $F(\tau, y, x)$ are of polynomial growth in $x\in \bbX$, uniformly in $t\in [0,T]$.
\item[(iv)] For each fixed $(\tau,y) \in [0, T)\times \bbX$ and $x\in\bbX$, $\mu^{\mathbf{u}}(t,x)$ and $\sigma^{\mathbf{u}}(t,x)$ are right-continuous in $t\in[0,T)$ and $\lim_{t'\ge t,(t',x')\rightarrow (t,x)}C^{\tau, y, \mathbf{u}}(t',x') = C^{\tau, y, \mathbf{u}}(t,x)$ for any $t\in[0,T)$.
\end{enumerate}
Denote the set of feasible strategies as $\mathbf{U}$.
\end{definition}

We impose the following assumption:
\begin{assumption}\label{as:ModelParameters}
Any $u\in \mathbb{U}$ is feasible.
\end{assumption}

\subsection{Intra-Personal Equilibrium}\label{subse:IntraPersonal}
Here and hereafter, $\hat{\mathbf{u}} \in \mathbf{U}$ denotes a given strategy and we examine whether it is an equilibrium strategy. For given $t\in[0,T)$, $\epsilon\in (0,T-t)$ and $\mathbf{a} \in \mathbf{U} $, define
\begin{align}\label{eq:PerturbatedPolicyFeedback}
{\mathbf{u}}_{t,\epsilon,\mathbf{a}}(s,y) := \begin{cases}
\mathbf{a}(s,y), &  s\in[t,t+\epsilon), y\in \bbX \\
\hat{\mathbf{u}}(s,y),&s\notin [t,t+\epsilon), y\in \bbX.
\end{cases}
\end{align}
Imagine that the agent at time $t$ chooses strategy $\mathbf{a}$ and is able to commit herself to this strategy in the period $[t,t+\epsilon)$. The agent, however, is unable to control her future selves beyond this small time period, namely in the period $[t+\epsilon,T)$ and believes that her future selves will take strategy $\hat{\mathbf{u}}$. Then, ${\mathbf{u}}_{t,\epsilon,\mathbf{a}}$ is the strategy that the agent at time $t$ expects herself to implement throughout the entire horizon. Note that $ {\mathbf{u}}_{t,\epsilon,\mathbf{a}}$ is feasible because both $\hat{\mathbf{u}}$ and $\mathbf{a}$ are feasible.

\begin{definition}[Intra-Personal Equilibrium]\label{de:EquilibriumFirstOrd}
$\hat{\mathbf{u}} \in \mathbf{U}$ is an {\em intra-personal equilibrium} if for any $x\in \bbX$, $t\in[0,T)$, and $\mathbf{a} \in \mathbf{U} $, we have
	\begin{align}\label{new requirement}
\limsup_{\epsilon\downarrow 0}\frac{ J(t,x;{\mathbf{u}}_{t,\epsilon, \mathbf{a}} )  -J(t,x;\hat{\mathbf{u}} ) }{\epsilon}\leq 0.
	\end{align}
\end{definition}

For each positive $\epsilon$, ${\mathbf{u}}_{t,\epsilon, \mathbf{a}}$ leads to a possibly different state process and thus to a different objective function value from those of $\hat{\mathbf{u}}$, so it is meaningful to compare the objective function values of ${\mathbf{u}}_{t,\epsilon, \mathbf{a}}$ and $\hat{\mathbf{u}}$ to examine whether the agent is willing to deviate from $\hat{\mathbf{u}}$ to $\mathbf{a}$ in the period of time $[t,t+\epsilon)$. Due to the continuous-time nature,  the length of the period,  $\epsilon$, during which the agent at $t$ exerts full self control, must be set to be infinitesimally small. Then, $J(t,x;{\mathbf{u}}_{t,\epsilon, \mathbf{a}} )$ and $J(t,x;\hat{\mathbf{u}} )$ become arbitrarily close to each other; so instead of evaluating their difference, we consider the {\em rate of increment} in the objective function value, i.e., the limit on the left-hand side of \eqref{new requirement}. Thus, under Definition \ref{de:EquilibriumFirstOrd}, a strategy $\hat{\mathbf{u}}$ is an intra-personal equilibrium if at any given time and state, the rate of increment in the objective value when the agent deviates from $\hat{\mathbf{u}}$ to any alternative strategy is nonpositive. As a result, the agent has little incentive to deviate from $\hat{\mathbf{u}}$.

\subsection{Sufficient and Necessary Condition}

We first introduce the generator of the controlled state process. Given $\mathbf{u} \in \mathbf{U}$ and interval $[a,b] \subseteq [0, T]$, consider $\xi$ that maps $(t,x)\in [a,b] \times \bbX$ to $ \xi(t,x)\in \bbR$.
Suppose $\xi \in \mathfrak{C}^{1,2}([a,b]\times \bbX)$, and denote by $\xi_t$, $\xi_x$, and $\xi_{xx}$ respectively
its first-order partial derivative in $t$, first-order partial derivative in $x$, and second-order partial derivative in $x$.
Define the following generator:
\begin{align}\label{eq:Generator}
\mathcal{A}^{\mathbf{u}}\xi(t,x)=\xi_t(t,x)+\xi_x(t,x)\mu^{\mathbf{u}}(t, x)+\frac{1}{2}\trace\left(\xi_{xx}(t,x)\tran \Upsilon^{\mathbf{u}}(t, x)\right),\notag \\
 t\in [a,b],x\in\bbX.
\end{align}

 For each fixed $(\tau,y)\in [0,T)\times \bbX$, denote
\begin{align}
&{f}^{\tau, y}(t,x):=\expect_{t,x}[F(\tau, y, X^{\hat{\mathbf{u}}}(T) ) ],\label{eq:fFunction} \\
 &{g}(t,x):=\expect_{t,x}[ X^{\hat{\mathbf{u}}}(T)  ],\; t\in[0,T],x\in\bbX.\label{eq:gFunction}
\end{align}
In addition, for fixed $(\tau,y) \in [0, T)\times \bbX$ and $s\in [0, T]$, denote
\begin{align}\label{eq:cFunction}
&{c}^{\tau,y, s}(t,x):=\expect_{t,x}[ C^{\tau, y, \hat{\mathbf{u}}}(s, X^{\hat{\mathbf{u}}}(s))],\; t\in[0,s],x\in\bbX.
\end{align}
In the following, $\mathcal{A}^{\mathbf{u}}f^{\tau,y}$ denotes the function that is obtained by applying the operator $\mathcal{A}^{\mathbf{u}}$ to $f^{\tau,y}(t,x)$ as a function of $(t,x)$ while fixing $(\tau,y)$.
Then, $\mathcal{A}^{\mathbf{u}}f^{t,x}(t,x)$ denotes the value of $\mathcal{A}^{\mathbf{u}}f^{\tau,y}$ at $(t,x)$ while $(\tau,y)$ is also set at $(t,x)$.
The above notations also apply to $C^{\tau,y,\mathbf{u}}$ and $c^{\tau,y,s}$.

To illustrate how to evaluate $J(t,x;{\mathbf{u}}_{t,\epsilon, \mathbf{a}} )  -J(t,x;\hat{\mathbf{u}} )$ and thus the rate of increment, let us consider the second term in the objective function (\ref{eq:ObjFun}). An informal calculation yields
\begin{align*}
  &\expect_{t,x}\left[F(t, x, X^{{\mathbf{u}}_{t,\epsilon,\mathbf{a}}}(T) )\right] - \expect_{t,x}\left[F(t, x, X^{\hat{\mathbf{u}}}(T) )\right]\\
  & =  \expect_{t,x}\left[\expect_{ t+\epsilon,X^{{\mathbf{u}}_{t,\epsilon,\mathbf{a}}}(t+\epsilon)}\left[F(t, x, X^{{\mathbf{u}}_{t,\epsilon,\mathbf{a}}}(T) )\right]\right] - \expect_{t,x}\left[F(t, x, X^{\hat{\mathbf{u}}}(T) )\right]\\
  &= \expect_{t,x}\left[f^{(t,x)}(t+\epsilon, X^{\mathbf{a}}(t+\epsilon) )\right] - f^{(t,x)}(t,x)\\
  & \approx \mathcal{A}^{\mathbf{a}}f^{t,x}(t,x)\epsilon,
\end{align*}
where the second equality holds because ${\mathbf{u}}_{t,\epsilon,\mathbf{a}}(s,\cdot)=\mathbf{a}(s,\cdot)$ for $s\in [t,t+\epsilon)$ and ${\mathbf{u}}_{t,\epsilon,\mathbf{a}}(s,\cdot) = \hat{\mathbf{u}}(s,\cdot)$ for $s\in [t+\epsilon,T)$ in addition to the definition of $f^{\tau,y}$ in \eqref{eq:fFunction}. The change of the other terms in the objective function when the agent deviates from $\hat{\mathbf{u}}$ to $\mathbf{a}$ in the period  $[t,t+\epsilon)$ can be evaluated similarly. As a result, we can derive the rate of increment in the objective value, namely the limit on the left-hand side of \eqref{new requirement}, which in turn enables us to derive a sufficient and necessary condition for $\hat{\mathbf{u}}$ to be an intra-personal equilibrium.

To formalize the above heuristic argument, we need to impose the following assumption:
\begin{assumption}\label{as:FirstOrderSmoothness}
	For any fixed $(\tau,y) \in [0, T)\times \bbX$ and $t\in[0,T)$, there exists $\tilde t\in (t,T]$ such that (i) ${f}^{\tau, y},g\in \mathfrak{\bar C}^{1,2}([t,\tilde t]\times\bbX)$; (ii) ${c}^{\tau, y, s}\in\mathfrak{C}^{1,2}([t,\tilde t\wedge s]\times \bbX)$ for each fixed $s\in (t,T]$ and
	$\frac{\partial^{j+\alpha} {c}^{\tau, y, s}(t',x') }{\partial t^j \partial x^\alpha}$ is of polynomial growth in $x'\in \bbX$, uniformly in $t'\in[t,\tilde t\wedge s]$ and $s\in (t,T]$, for any $\alpha$ with $|\alpha|\le 2-2j$ and $j=0,1$; and (iii) $G(\tau,y,z)$ is continuously differentiable in $z$, with the partial derivative denoted as $G_z(\tau,y,z)$.
\end{assumption}

\begin{theorem}\label{Theorem:first order derivative}
	Suppose Assumptions \ref{as:ModelParameters} and \ref{as:FirstOrderSmoothness} hold. Then, for any $(t,x)\in [0,T)\times \bbX$ and $\mathbf{a} \in \mathbf{U}$,  we have
	\begin{align}\label{first order expansion}
	&\lim_{\epsilon\downarrow 0}\frac{J(t,x;{\mathbf{u}}_{t,\epsilon,\mathbf{a}} )-J(t,x;\hat{\mathbf{u}} )}{\epsilon}= \Gamma^{t,x,\hat{\mathbf{u}}}(t,x; \mathbf{a}),
	\end{align}
where for any $(\tau,y)\in [0,T)\times \bbX$,
	\begin{align}\label{eq:FirstOrderDerivative}
	\Gamma^{\tau, y,\hat{\mathbf{u}}}(t,x; \mathbf{a}):&=C^{\tau, y, \mathbf{a}}(t,x )-{C}^{\tau, y, \hat{\mathbf{u}}}(t,x)
	+\int_t^T \mathcal{A}^{\mathbf{a}} {c}^{\tau,y,s}(t, x )ds\notag \\
&+ \mathcal{A}^{\mathbf{a} }  f^{\tau,y}(t,x) +G_z(\tau,y,  g(t,x))\mathcal{A}^{\mathbf{a} }  g(t,x).
	\end{align}
	Moreover, $\Gamma^{\tau,y,\hat{\mathbf{u}}}(t,x; \mathbf{a})=\Gamma^{\tau,y,\hat{\mathbf{u}}}(t,x; \tilde{\mathbf{a}})$ for any  $\mathbf{a}, \tilde{\mathbf{a}} \in  \mathbf{U}$ with $\mathbf{a}(t,x)=\tilde{\mathbf{a}}(t,x)$ and $\Gamma^{\tau,y,\hat{\mathbf{u}}}(t,x; \mathbf{a})=0$ if $\mathbf{a}(t,x)=\hat{\mathbf{u}}(t,x)$.
	Consequently, $\hat{\mathbf{u}}$ is  an intra-personal  equilibrium if and only if
	\begin{align}\label{weak equilibrium condition}
\Gamma^{t,x,\hat{\mathbf{u}}}(t,x; u)\leq 0, \; \forall u\in \mathbb{U}, x\in \bbX, t\in [0, T).
	\end{align}
\end{theorem}

Theorem \ref{Theorem:first order derivative} presents a sufficient and necessary condition \eqref{weak equilibrium condition} for an intra-personal equilibrium $\hat{\mathbf{u}}$. Because $\Gamma^{\tau,y,\hat{\mathbf{u}}}(\tau,y; \hat{\mathbf{u}}(t,x))=0$, we have
\begin{align*}
  \Gamma^{\tau, y,\hat{\mathbf{u}}}(t,x; \mathbf{a}) = \Pi^{\tau,y}(t,x;\mathbf{a})-\Pi^{\tau,y}(t,x;\hat{\mathbf{u}}),
\end{align*}
where
\begin{align}
  \Pi^{\tau,y}(t,x;\mathbf{a}):&=C^{\tau, y, \mathbf{a}}(t,x ) + \int_t^T \mathcal{A}^{\mathbf{a}} {c}^{\tau,y,s}(t, x )ds+ \mathcal{A}^{\mathbf{a} }  f^{\tau,y}(t,x)\notag\\
  &+ G_z(\tau,y,  g(t,x))\mathcal{A}^{\mathbf{a} }  g(t,x).\label{eq:Pifun}
\end{align}
As a result, condition \eqref{weak equilibrium condition} is equivalent to
	\begin{align}\label{eq:EquiConditionForm2}
\max_{u\in \mathbb{U}}\Gamma^{t,x,\hat{\mathbf{u}}}(t,x; u)=0, x\in \bbX, t\in [0, T)
\end{align}
or
	\begin{align}\label{eq:EquiConditionForm3}
\hat{\mathbf{u}}(t,x)\in\mathrm{arg}\max_{u\in \mathbb{U}}\Pi^{t,x}(t,x;u), x\in \bbX, t\in [0, T).
\end{align}
This can be regarded as a time-inconsistent version of the verification theorem
in (classical) stochastic control.

The proof of Theorem \ref{Theorem:first order derivative} can be found in \citet{Bjork2017:TimeInconsistent} and \citet{HeJiang2019:OnEquilibriumStrategies}. Assumption \ref{as:ModelParameters} is easy to verify because it involves only the model parameters, i.e., $\mu$, $\sigma$, $C$, $F$, and $G$. Assumption \ref{as:FirstOrderSmoothness} imposes some regularity conditions on $\hat{\mathbf{u}}$, which usually requires $\hat{\mathbf{u}}$ to be smooth to a certain degree; see \citet{HeJiang2019:OnEquilibriumStrategies} for a sufficient condition for this assumption. As a result, the sufficient and necessary condition \eqref{weak equilibrium condition} cannot tell us whether there exists any intra-personal equilibrium among the strategies that do not satisfy Assumption \ref{as:FirstOrderSmoothness}. This  condition, however, is still very useful for us to find intra-personal equilibria for specific problems. Indeed, in most time-inconsistent problems in the literature, intra-personal equilibrium can be found and verified using \eqref{weak equilibrium condition}; see Section \ref{sec:Applications}.

\subsection{Extended HJB}\label{subse:ExtendedHJB}


Define the {\em continuation value} of a strategy $\hat{\mathbf{u}}$, denoted as $V^{\hat{\mathbf{u}}}(t,x),(t,x)\in [0,T]\times \bbX$, to be the objective value over time and state under this strategy, i.e.,
\begin{align}
  V^{\hat{\mathbf{u}}}(t,x):&=J(t,x;\hat{\mathbf{u}} ) = H^{t,x}(t,x) + G\big(t, x,g(t,x)  \big),\label{eq:ContinuationValue}
\end{align}
where
\begin{align}
  H^{\tau,y}(t,x):&=\expect_{t,x}\left[\int_t^T C^{\tau, y, \hat{\mathbf{u}}}(s, X^{\hat{\mathbf{u}}}(s))ds + F(\tau,y,X^{\hat{\mathbf{u}}}(T))\right]\notag\\
  & = \int_t^Tc^{\tau,y,s}(t,x)ds + f^{\tau,y}(t,x).\label{eq:Hfunction}
\end{align}
Assuming certain regularity conditions and applying the operator ${\cal A}^{u}$ to $V^{\hat{\mathbf{u}}}(t,x)$, we derive
\begin{align*}
  {\cal A}^{u} V^{\hat{\mathbf{u}}}(t,x) &= -C^{t,x,\hat{\mathbf{u}}}(t,x) + \int_t^T{\cal A}^{u} c^{t,x,s}(t,x)ds + {\cal A}^{u}f^{t,x}(t,x)\\
  & + G_{z}\big(t,x,g(t,x)\big){\cal A}^{u} g(t,x) + {\cal A}^{u}_{\tau, y} H^{t,x}(t,x)+{\cal A}^{u}_{\tau, y}G(t,x,g(t,x))\\
  & + \trace\left(\left(H^{t,x}_{xy}(t,x)+ G_{zy}(t,x,g(t,x))\tran g_x(t,x)\right)\tran \Upsilon^{u}(t, x) \right)\\
  & +\frac{1}{2} G_{zz}\big(t, x,g(t,x)  \big)\trace\left(g_x(t,x)g_x(t,x)\tran \Upsilon^{u}(t, x)\right)
\end{align*}
where $H^{\tau,y}_{xy}(t,x)$ denotes the cross partial derivative of $H^{\tau,y}(t,x)$ in $x$ and $y$, $G_{zy}(\tau,y,z)$ the cross partial derivative of $G(\tau,y,z)$ in $z$ and $y$, and $G_{zz}(\tau,y,z)$ the second-order derivative of $G(\tau,y,z)$ in $z$. For each {\em fixed} $(t,x)$, ${\cal A}^{u}_{\tau, y}H^{\tau,y}(t,x)$ denotes the generator of ${\cal A}^{u}$ applied to $H^{\tau,y}(t,x)$ as a function of $(\tau,y)$, i.e., ${\cal A}^{u}_{\tau, y}H^{\tau,y}(t,x):={\cal A}^u\ell (\tau,y)$, where $\ell(\tau,y):=H^{\tau,y}(t,x),(\tau,y)\in [0,T)\times \bbX$, and ${\cal A}^{u}_{\tau, y}G(\tau,y,g(t,x))$ is defined similarly.

Now, suppose $\hat {\mathbf{u}}$ is an intra-personal equilibrium. Recalling \eqref{eq:FirstOrderDerivative} and the sufficient and necessary condition \eqref{eq:EquiConditionForm2}, we derive the following equation satisfied by the continuation value of an intra-personal equilibrium $\hat{\mathbf{u}}$:
\begin{align}
  &\max_{u\in\mathbb{U}} \Big[{\cal A}^{u} V^{\hat{\mathbf{u}}}(t,x) + C^{t,x,u}(t,x)-\big({\cal A}^{u}_{\tau, y} H^{t,x}(t,x)+{\cal A}^{u}_{\tau, y}G(t,x,g(t,x))\big)\notag \\
  &-\trace\left(\left(H^{t,x}_{xy}(tx)+ G_{zy}(t,x,g(t,x))\tran g_x(t,x)\right)\tran \Upsilon^{u}(t, x) \right)\notag \\
 & -\frac{1}{2} G_{zz}\big(t, x,g(t,x)  \big)\trace\left(g_x(t,x)g_x(t,x)\tran \Upsilon^{u}(t, x)\right)\Big]=0, (t,x)\in [0,T)\times \bbX,\notag \\
  &V^{\hat{\mathbf{u}}}(T,x)  =  F(T,x,x) + G(T,x,x),\; x\in \bbX.\label{eq:EHJBEqMain}
\end{align}
By \eqref{eq:Hfunction}, the definitions of $c^{\tau,y,s}(t,x)$ and $f^{\tau,y}(t,x)$, and the Feymann-Kac formula, we derive the following equation for $H^{\tau,y}(t,x)$:
\begin{align}
  &{\cal A}^{\hat{\mathbf{u}}}H^{\tau,y}(t,x) + {C}^{\tau, y, \hat{\mathbf{u}}}(t,x) = 0, \quad (t,x)\in [0,T)\times \bbX,(\tau,y)\in [0,T)\times \bbX,\notag\\
  & H^{\tau,y}(T,x) = F(\tau,y,x),\quad x\in \bbX,(\tau,y)\in [0,T)\times \bbX.\label{eq:EHJBEqExpt}
\end{align}
Similarly, we derive the following equation for $g$:
\begin{align}
  &{\cal A}^{\hat{\mathbf{u}}}g(t,x) = 0, \quad (t,x)\in [0,T)\times \bbX,\notag\\
  & g(T,x) = x,\quad x\in \bbX.\label{eq:EHJBEqMean}
\end{align}

Some remarks are in order. First, instead of a single equation for the value function of a time-consistent problem, the intra-personal equilibrium and its continuation value satisfy a system of equations \eqref{eq:EHJBEqMain}--\eqref{eq:EHJBEqMean}, which is referred to  as the {\em extended HJB equation} by \citet{Bjork2017:TimeInconsistent}.

Second, compared to the HJB equation for a time-consistent problem, which takes the form $\max_{u\in\mathbb{U}} \big[{\cal A}^{u} V^{\hat{\mathbf{u}}}(t,x) + C^{u}(t,x)\big]=0$, equation \eqref{eq:EHJBEqMain} has three additional terms in the first, second, and third lines of the equation, respectively. Here and hereafter, when $C^{\tau,y,u}(t,x)$ does not depend on $(\tau,y)$, we simply drop the superscript $(\tau,y)$. Similar notations apply to $H^{\tau,y}(t,x)$ and to the case when there is no dependence on $y$. Now, recall that for the objective function \eqref{eq:ObjFun}, time inconsistency arises from (i) the dependence of $C$, $F$, and $G$ on $(t,x)$ and (ii) the nonlinear dependence of $G\big(t, x, \expect_{t,x} [  X^{{\mathbf{u}}}(T) ]\big)$ on $\expect_{t,x} [  X^{{\mathbf{u}}}(T) ]$. If source (i) of time inconsistency is absent, the first and second additional terms in \eqref{eq:EHJBEqMain} will vanish. If source (ii) of time inconsistency is absent, the third additional term in \eqref{eq:EHJBEqMain} will disappear. In particular, without time inconsistency, the extended HJB equation \eqref{eq:EHJBEqMain} reduces to the classical HJB equation.

Third, consider the case in which $C$, $F$, and $G$ do not depend on $x$ and $G\big(t, x, \expect_{t,x} [  X^{{\mathbf{u}}}(T) ]\big)$ is linear in $\expect_{t,x} [  X^{{\mathbf{u}}}(T) ]$. In this case, the second and third lines of \eqref{eq:EHJBEqMain} vanish and we can assume $G\equiv 0$ without loss of generality because $G$ can be combined with $F$. As a result, the extended HJB equation \eqref{eq:EHJBEqMain} specializes to
\begin{align}
  &\max_{u\in\mathbb{U}} \Big[{\cal A}^{u} V^{\hat{\mathbf{u}}}(t,x) + C^{t,u}(t,x)\Big]=h^{t}(t,x), (t,x)\in [0,T)\times \bbX,\notag \\
  &V^{\hat{\mathbf{u}}}(T,x)  =  F(T,x),\; x\in \bbX\label{eq:EHJBEqSimplified}
\end{align}
where $h^{\tau}(t,x):=H^{\tau}_\tau(t,x)$ (with the subscript $\tau$ denoting the partial derivative with respect to $\tau$) and thus satisfies
\begin{align}
  &{\cal A}^{\hat{\mathbf{u}}}h^{\tau}(t,x) + {C}^{\tau,\hat{\mathbf{u}}}_\tau(t,x) = 0, \quad (t,x)\in [0,T)\times \bbX,\tau \in[0,T),\notag\\
  & h^{\tau}(T,x) = F_\tau(\tau,x),\quad x\in \bbX,\tau\in [0,T).\label{eq:EHJBEqExptDer}
\end{align}

\section{Discussions}\label{sec:Discussions}
\subsection{Intra-Personal Equilibria with Fixed Initial Data}

Consider an agent at time 0 with a fixed state $x_0$ who correctly anticipates that her self at each future time $t$ faces the problem \eqref{eq:ControlProblem} and who has no control of future selves at any time. A strategy $\hat{\mathbf{u}}$ can be consistently implemented by the agent throughout the entire horizon $[0,T]$ if the agent has no incentive to deviate from it at any time {\em along the state path}.  Actions that the agent might be taking were she not on the state path are irrelevant. To be more precise, for any fixed initial data $(0,x_0)$, we define $\hat{\mathbf{u}}$ to be an {\em intra-personal equilibrium starting from $(0,x_0)$} if \eqref{new requirement} holds for any $\mathbf{a}\in \mathbf{U}$, $t\in [0,T)$, and $x\in \mathbb{X}^{0,x_0,\hat{\mathbf{u}}}_t$, where $\mathbb{X}^{0,x_0,\hat{\mathbf{u}}}_t$ denotes the set of all possible states at time $t$ along the state path starting from $x_0$ at the initial time and under the strategy $\hat{\mathbf{u}}$. 

It is evident that the intra-personal equilibrium defined in Definition \ref{de:EquilibriumFirstOrd} is {\em universal} in that it is an  equilibrium starting from {\em any} initial data $(0,x_0)$. On the other hand, starting from a fixed state $x_0$ at time 0, the state process in the future might not be able to visit the whole state space; so an  equilibrium starting from $(0,x_0)$ is not necessarily universal, i.e., it is not necessarily an  equilibrium when the agent starts from other initial data. For example, \citet{HeEtal2019:MedianMaximization} consider a continuous-time portfolio selection problem in which an agent maximizes the median of her terminal wealth. With a fixed initial wealth of the agent, the authors derive a set of intra-personal equilibrium strategies starting from this particular initial wealth level. They show that these strategies are no longer equilibria if the agent starts from some other initial wealth levels, and in particular not universal equilibria in the sense of Definition \ref{de:EquilibriumFirstOrd}.


The first study of intra-personal equilibria starting from a fixed initial data dates back to \cite{PelegYaari1973:OnExistenceConsistent}. In a discrete-time setting, the authors propose that a strategy $(s^*_0,s^*_1,\dots)$, where $s^*_t$ stands for the agent's closed-loop strategy at time $t$, is an equilibrium strategy if for any $t$, $(s_0^*,\dots, s_{t-1}^*,s_t,s_{t+1}^*,\dots )$ is dominated by $(s_0^*,\dots, s_{t-1}^*,s_t^*,s_{t+1}^*,\dots )$ for any $s_t$. They argue that the above definition is more desirable than the following one, which is based on a model in \cite{Pollak1968:ConsistentPlanning}: $(s^*_0,s^*_1,\dots)$ is an equilibrium strategy if for any time $t$, $(s_0,\dots, s_{t-1},s_t,s_{t+1}^*,\dots )$ is dominated by $(s_0,\dots, s_{t-1},s_t^*,s_{t+1}^*,\dots )$ for {\em any} $(s_0,\dots, s_t)$. It is clear that the equilibrium strategies considered by \cite{PelegYaari1973:OnExistenceConsistent} are the ones starting from a fixed initial data while those studied by \cite{Pollak1968:ConsistentPlanning} are universal. Recently, \citet{HeJiang2019:OnEquilibriumStrategies}, \citet{HanWong2020:TimeInconsistency}, and \citet{HernandezPossami2020:MyMyself} also consider intra-personal equilibria with fixed initial data. Moreover, \citet{HeJiang2019:OnEquilibriumStrategies} propose a formal definition of $\mathbb{X}^{0,x_0,\hat{\mathbf{u}}}_t$, calling it the set of reachable states.

Finally, let us comment that the sufficient and necessary condition in Theorem \ref{Theorem:first order derivative} is still valid for intra-personal equilibria starting from fixed initial data $(0,x_0)$, provided that we replace $\bbX$ in this condition with the set of reachable states $\mathbb{X}^{0,x_0,\hat{\mathbf{u}}}_t$; see  \citet{HeJiang2019:OnEquilibriumStrategies} for details. The extended HJB equation in Section \ref{subse:ExtendedHJB} can be revised and applied similarly.

\subsection{Set of Alternative Strategies}
In Definition \ref{de:EquilibriumFirstOrd}, the set of strategies that the agent can choose at time $t$ to implement for the period $[t,t+\epsilon)$, denoted as $\mathbf{D}$, is set to be the entire set of feasible strategies $\mathbf{U}$. This definition is used in \citet{Bjork2017:TimeInconsistent}, \cite{EkelandPirvu2008:InvestConsptnWithoutCommit}, and \cite{EkelandEtal2012:TimeConsistent}. In some other works, however, $\mathbf{D}$ is set to be the set of constant strategies $\mathbb{U}$; see for instance \cite{EkelandLazrak2006:BeingSeriousAboutNonCommitment,EkelandILazrakA:07ti,EkelandLazrak2010:GoldenRule}, \cite{BjordTMurgoci:08u}, and \cite{BasakSChabakauri2010:DynamicMeanVariance}. \citet{HeJiang2019:OnEquilibriumStrategies} show that the choice of $\mathbf{D}$ is irrelevant as long as it at least contains $\mathbb{U}$. Indeed, this can be seen from the observation in Theorem \ref{Theorem:first order derivative} that $\Gamma^{\tau,y,\hat{\mathbf{u}}}(t,x; \mathbf{a})=\Gamma^{\tau,y,\hat{\mathbf{u}}}(t,x; \mathbf{a}(t,x))$ for any  $\mathbf{a} \in \mathbf{U}$. \citet{HeJiang2019:OnEquilibriumStrategies} also show that for strong intra-personal equilibrium, which will be introduced momentarily, the choice of $\mathbf{D}$ is relevant.

\subsection{Regular and Strong Intra-Personal Equilibrium}
As noted in Remark 3.5 of \cite{Bjork2017:TimeInconsistent}, condition \eqref{new requirement} does not necessarily imply that $J(t,x;{\mathbf{u}}_{t,\epsilon, \mathbf{a}} )$ is less than or equal to $J(t,x;\hat{\mathbf{u}} )$ however small $\epsilon>0$ might be and thus disincentivizes   the agent from deviating from $\hat{\mathbf{u}}$. For example, if $J(t,x;{\mathbf{u}}_{t,\epsilon, \mathbf{a}} )-J(t,x;\hat{\mathbf{u}} )=\epsilon^2$, then \eqref{new requirement} holds, but the agent can achieve a strictly larger objective value if she deviates from $\hat{\mathbf{u}}$ to $\mathbf{a}$ and thus is willing to do so. 

To address the above issue, \cite{HuangZhou2018:StrongWeakEquilibria} and \citet{HeJiang2019:OnEquilibriumStrategies} propose the notion of {\em strong intra-personal equilibrium}:
\begin{definition}[Strong Intra-personal Equilibrium]\label{de:EquilibriumExtd}
$\hat{\mathbf{u}} \in \mathbf{U}$ is a {\em strong intra-personal equilibrium strategy}  if for any $x\in \bbX$, $t\in[0,T)$, and $\mathbf{a} \in \mathbf{D}$, there exists $\epsilon_0 \in (0,T-t)$ such that
\begin{align}\label{new requirementExtd}
J(t,x;{\mathbf{u}}_{t,\epsilon,\mathbf{a}} )  -J(t,x;\hat{\mathbf{u}} ) \leq 0,\quad \forall \epsilon\in(0,\epsilon_0].
\end{align}
\end{definition}
It is straightforward to see that a strong intra-personal equilibrium implies the one in Definition \ref{de:EquilibriumFirstOrd}, which we refer to as a {\em weak intra-personal equilibrium} in this subsection.

\cite{HuangZhou2018:StrongWeakEquilibria} consider a stochastic control problem in which an agent can control the generator of a time-homogeneous, continuous-time, finite-state Markov chain at each time to maximize expected running reward in an infinite time horizon. Assuming that at each time the agent can implement a time-homogeneous strategy only, the authors provide a characterization of a strong intra-personal equilibrium and prove its existence under certain conditions.

\citet{HeJiang2019:OnEquilibriumStrategies} follow the framework in \eqref{sec:ExHJB} and derive two necessary conditions for a strategy to be strong intra-personal equilibrium. Using these conditions, the authors show that strong intra-personal equilibrium does not exist for the portfolio selection and consumption problems studied in \citet{EkelandPirvu2008:InvestConsptnWithoutCommit}, \citet{BasakSChabakauri2010:DynamicMeanVariance}, and \citet{BjorkEtal2011:MeanVariancewithStateDependentRiskAversion}. Motivated by this non-existence result, the authors propose the so-called {\em regular intra-personal equilibrium} and show that it exists for the above three problems and is stronger than the weak intra-personal equilibrium and weaker than the strong intra-personal equilibrium in general.

\subsection{Existence and Uniqueness}
In most studies on time-inconsistent problems in the literature, a {\em closed-form} strategy is constructed and verified to satisfy the sufficient and necessary condition \eqref{weak equilibrium condition} or the extended HJB equation \eqref{eq:EHJBEqMain}--\eqref{eq:EHJBEqMean}. The existence of intra-personal equilibrium in general is difficult to prove because it essentially relies on  a fixed point argument: For each guess of intra-personal equilibrium $\hat{\mathbf{u}}$, we first calculate $\Gamma^{\tau,y,\hat{\mathbf{u}}}$ in \eqref{weak equilibrium condition} and $H^{\tau,y}(t,x)$ and $g$ in \eqref{eq:EHJBEqExpt} and \eqref{eq:EHJBEqMean}, respectively, and then derive an updated intra-personal equilibrium, denoted as $\mathbb{T} \hat{\mathbf{u}}$, from the condition \eqref{weak equilibrium condition} or from the equation \eqref{eq:EHJBEqMain}. The existence of an intra-personal equilibrium then  boils down to the existence of the fixed point of $\mathbb{T}$. The mapping $\mathbb{T}$ is highly nonlinear; so the existence of its fixed point is hard to establish. Additional difficulty is caused by the regularity conditions that we need to pose on $\hat{\mathbf{u}}$ to validate the sufficient and necessary condition \eqref{weak equilibrium condition} or the extended HJB equation \eqref{eq:EHJBEqMain}--\eqref{eq:EHJBEqMean}.

We are only aware of very few works on the existence of intra-personal equilibria in continuous time. \citet{Yong2012:TimeInconsistent} proposes an alternative approach to defining the strategy of a sophisticated agent, which will be discussed in detail in Section \ref{sec:Discretization}. Assuming $G\equiv 0$, $C$ and $F$ to be independent of $x$ in the objective function \eqref{eq:ObjFun}, and $\sigma(t,x,u)$ in the controlled diffusion process \eqref{eq:ControlProblem} to be independent of control $u$ and nondegenerate, \citet{Yong2012:TimeInconsistent} proves the existence of the sophisticated agent's strategy, which is used to imply the existence of an intra-personal equilibrium under Definition \ref{de:EquilibriumFirstOrd}. \citet{WeiEtal2017:TimeInconsistent} and \citet{WangYong2019:TimeInconsistent} extend the result of \citet{Yong2012:TimeInconsistent} by generalizing the objective function; however  for the existence of intra-personal equilibria, they need to assume the volatility $\sigma$ to be independent of control and nondegenerate. \citet{HernandezPossami2020:MyMyself} study intra-personal equilibria in a non-Markovian setting, where they consider a non-Markovian version of the objective function in \citet{Yong2012:TimeInconsistent} and assume the drift $\mu$ of the controlled process to be in the range of the volatility matrix at each time. The authors prove the existence of intra-personal equilibria when the volatility $\sigma$ is independent of control.

Intra-personal equilibria  can be non-unique; see \citet{EkelandLazrak2010:GoldenRule}, \citet{CaoWerning2016:DynamicSavingsDisagreements}, and \citet{HeEtal2019:MedianMaximization}. For some problems, however, uniqueness has been established in the literature. Indeed, \citet{Yong2012:TimeInconsistent}, \citet{WeiEtal2017:TimeInconsistent}, \citet{WangYong2019:TimeInconsistent}, and \citet{HernandezPossami2020:MyMyself} prove the uniqueness in various settings with the common assumption that the volatility $\sigma$ is independent of control.

%

\subsection{Non-Markovian Strategies}
In most studies on time-inconsistent problems, where the controlled state processes are Markovian, the search for intra-personal equilibrium is restricted to the set of Markovian strategies, i.e., strategies that are functions of time $t$ and the {\it current} state value $x$. Motivated by some practical problems  such as rough volatility models and principle-agent problems, \citet{HanWong2020:TimeInconsistency} and \citet{HernandezPossami2020:MyMyself} define and search intra-personal equilibria in the class of non-Markovian or path-dependent strategies, i.e., ones that depend on time $t$ and the whole path of the controlled state up to $t$.

\section{Closed-Loop versus Open-Loop Intra-Personal Equilibria}\label{sec:OpenLoop}

A {\em closed-loop} or {\it feedback} control strategy is a function $\mathbf{u}$ that maps time $t$ and the controlled state path $(x_s)_{s\le t}$ up to $t$ to the space of actions. As a result, the action taken by an agent under such a strategy is $\mathbf{u}(t, (x_s)_{s\le t})$. An {\em open-loop} control is a collection of actions over time and state of the nature, $(u(t,\omega))_{t\ge 0}$, where $u(t,\omega)$ is the action to be taken at time $t$ and in scenario $\omega$, {\em regardless} of the state path $(x_s)_{s\le t}$. For classical time-consistent control problems and under some technical assumptions, the state-control paths under the optimal open-loop control and under the optimal closed-loop control strategy are the same if the controlled system starts from the same initial time and state; see for instance \citet{YongJZhouXY:99sc}.

In Section \ref{sec:ExHJB}, intra-personal equilibrium is defined for closed-loop control strategies, which is also the approach taken by most studies on time-inconsistent problems in the literature. In some other works, intra-personal equilibrium is defined for open-loop controls; see for instance \citet{HuEtal2012:TimeInconsistentStochasticLQ}, \citet{HuEtal2017:TimeInconsisten}, \citet{LiEtal2019:EquilibriumStrategies}, and \citet{HuEtal2020:ConsistentInvestmentRDU}.

Formally, under the same probabilistic framework in Section \ref{subse:framework}, we represent an open-loop strategy by a progressively measurable process $(u(t))_{t\ge 0}$ that takes values in $\mathbb{U}$. The controlled state process $X^{u}$ takes the form
\begin{align*}
dX^{u}(s)=\mu(s, X^{u}(s),u(s) )ds  +\sigma(s,X^{u}(s), u(s) )dW(s),\; s\in[t,T];\; X^{u}(t)=x.
\end{align*}
Denote by ${\cal U}$ the set of feasible open-loop controls, i.e., the set of progressively measurable processes on $[0,T]$ satisfying certain integrability conditions. At time $t$ with state $x$, suppose the agent's objective is to maximize $J(t,x;u(\cdot))$ by choosing $u(\cdot)\in {\cal U}$. Given $\hat u(\cdot)\in {\cal U}$, for any $t\in [0,T)$, $x\in \bbX$, $\epsilon \in (0,T-t)$, and $a(\cdot)\in {\cal U}$, define
\begin{align}\label{eq:PerturbatedPolicyOpenLoop}
u_{t,\epsilon,a}(s) := \begin{cases}
a(s), &  s\in[t,t+\epsilon) \\
\hat u(s),&s\notin [t,t+\epsilon).
\end{cases}
\end{align}
Suppose that at time $t$ with state $x$, the agent chooses an open-loop control $a(\cdot)$, but is only able to implement it in the period $[t,t+\epsilon)$. Anticipating that her future selves will take the given control $\hat u(\cdot)$, the agent expects herself to follow $u_{t,\epsilon,a}$ in the period $[t,T]$. 
\begin{definition}[Open-Loop Intra-Personal Equilibrium]\label{de:EquilibriumOpenLoop}
$\hat u(\cdot)\in {\cal U}$ is an {\em open-loop intra-personal equilibrium} if for any $x\in \bbX$, $t\in[0,T)$, and $a \in {\cal U}$ that is constant in a small period after $t$, we have
	\begin{align}\label{eq:OpenLoopEquilibriumCond}
\limsup_{\epsilon\downarrow 0}\frac{ J(t,x;u_{t,\epsilon,a}(\cdot) )  -J(t,x;\hat u(\cdot) ) }{\epsilon}\leq 0.
	\end{align}
\end{definition}

The above is analogous to  the definition of an intra-personal equilibrium for closed-loop strategies. However, there is a subtle yet crucial difference between the two definitions. For the one for open-loop controls, the perturbed control $u_{t,\epsilon,a}(s)$ defined by (\ref{eq:PerturbatedPolicyOpenLoop}) and the original one $\hat u$ are identical on $[t+\epsilon,T]$ as two stochastic processes. In other words, the perturbation in the small time period $[t,t+\epsilon)$ will not affect the control process beyond this period. This is not the case for the closed-loop counterpart, because the perturbation (\ref{eq:PerturbatedPolicyFeedback}) on $[t,t+\epsilon)$ changes the 
 control in the period,
 which will alter the state process in $[t,t+\epsilon)$ and in particular the
state at time $t+\epsilon$. This in turn will change the control {\it process} on $[t+\epsilon,T]$ upon substituting  the state process  into the feedback strategy.

To characterize open-loop intra-personal equilibria, we only need to compute the limit on the left-hand side of \eqref{eq:OpenLoopEquilibriumCond}. This limit can be evaluated by applying the spike variation technique that is used to derive Pontryagin's maximum principle for time-consistent control problems in continuous time \citep{YongJZhouXY:99sc}. As a result, open-loop intra-personal equilibrium can be characterized by a flow of forward-backward stochastic differential equations (SDEs); see \citet{HuEtal2012:TimeInconsistentStochasticLQ} for more details.
In contrast, the spike variation technique no longer works for closed-loop equilibria
because the perturbed control process is different from the original one beyond the small time period for perturbation, as discussed above.


This discussion suggests  that closed-loop and open-loop  equilibria are likely different. This is confirmed by \citet{HuEtal2012:TimeInconsistentStochasticLQ}. The authors consider a mean-variance portfolio selection problem, where an agent decides the dollar amount invested in a stock at each time, and derive an open-loop equilibrium; see Section 5.4.1 therein. They then compare this equilibrium with the closed-loop  equilibrium derived by \citet{BjorkEtal2011:MeanVariancewithStateDependentRiskAversion} for the same portfolio selection problem, and find that the state-control path under these two equilibria are different.

It can be argued that closed-loop strategies are preferred to the open-loop ones for three reasons. First, in many problems, agents' actions naturally depend on some state variables. For example, in a consumption problem, an agent's consumption at any time is more likely to depend directly on her wealth at that time. If her wealth suddenly increases, she would  probably consume more.

Second, closed-loop intra-personal equilibrium is invariant to the choice of control variables while open-loop intra-personal equilibrium might not. For example, in a portfolio selection problem where an agent decides the allocation of her wealth between a risk-free asset and a risky stock, the decision variable can be set to be the dollar amount invested in the stock or the percentage of wealth invested in the stock. Suppose $\hat {\mathbf{u}}$ is a closed-loop intra-personal equilibrium representing the percentage of wealth invested in the stock. Then, we have
\begin{align}\label{eq:IntraPEquiClosedLoopCond}
\limsup_{\epsilon\downarrow 0}\frac{ J(t,x;{\mathbf{u}}_{t,\epsilon, \mathbf{a}} )  -J(t,x;\hat{\mathbf{u}} ) }{\epsilon}\leq 0.
\end{align}
for all $t\in [0,T)$, $x\in \bbX$, and $\mathbf{a}\in \mathbf{U}$, where the state variable $x$ represents the agent's wealth. Now, suppose we represent the agent's decision by the dollar amount invested in the risky stock, and denote a control strategy as $\boldsymbol{\pi}$. Then, the agent's objective function is $\tilde J(t,x;\boldsymbol{\pi}) = J(t,x;\mathbf{u})$ with $\mathbf{u}(s,y) = \boldsymbol{\pi}(s,y)/y$. Condition \eqref{eq:IntraPEquiClosedLoopCond} implies that
\begin{align*}
  \limsup_{\epsilon\downarrow 0}\frac{ \tilde J(t,x;\boldsymbol{\pi}_{t,\epsilon, \tilde{\mathbf{a}}} )  -\tilde J(t,x;\hat{\boldsymbol{\pi}} ) }{\epsilon}\leq 0,
\end{align*}
for any $t\in [0,T)$, $x\in \bbX$, and strategy $\tilde{\mathbf{a}}$ that represents the dollar amount invested in the stock, where $\hat{\boldsymbol{\pi}}(s,y):=y\hat{\mathbf{u}}(s,y)$ and $\boldsymbol{\pi}_{t,\epsilon, \tilde{\mathbf{a}}}$ is defined similarly to ${\mathbf{u}}_{t,\epsilon, \mathbf{a}}$. Thus, $\hat{\boldsymbol{\pi}}$, which is the dollar amount investment strategy corresponding to the percentage investment strategy $\hat {\mathbf{u}}$, is also an intra-personal equilibrium. By contrast, for the mean-variance portfolio selection problem studied by \citet{HuEtal2012:TimeInconsistentStochasticLQ}, where the agent's decision is the dollar amount invested in the stock, the open-loop intra-personal equilibrium yields a different control-state path from the one yielded by its closed-loop counterpart  derived by \citet{BjorkEtal2011:MeanVariancewithStateDependentRiskAversion}. If we change the agent's decision variable to the percentage of wealth invested in the stock, the open-loop intra-personal equilibrium and the closed-loop intra-personal equilibrium in \citet{BjorkEtal2011:MeanVariancewithStateDependentRiskAversion} yield the same control-state path. This implies that open-loop equilibria depend on  the choice of control variables.

Third, open-loop intra-personal equilibrium may not be well-posed for some problems. Consider the discrete-time version of the consumption problem studied in \citet{Strotz1955:MyopiaInconsistency}: An agent decides the amount of consumption $C_t$ at each time $t=0,1,\dots, T$ with the total budget $x_0$, i.e., $\sum_{t=0}^T C_t =x_0$. For this problem, any consumption plan $(\hat C_t)_{t\ge 0}$ is an open-loop intra-personal equilibrium. Indeed, at each time $t$, anticipating her future selves will consume $\hat C_s,s=t+1,\dots, T$, the only amount of consumption $C_t$ that the agent can choose at time $t$ is $\hat C_t$ due to the budget constraint $(\sum_{s=0}^{t-1}\hat C_s) + C_t + (\sum_{s=t+1}^T\hat C_s)=x_0$. This leads to a trivial definition of intra-personal equilibrium. The above issue can be rectified  if we use closed-loop strategies. To see this, we set $x_t$ to be the agent's remaining budget at time $t$ before the consumption at that time. For closed-loop intra-personal equilibrium, we consider a mapping from time $t$ and the remaining budget $x_t$ to the consumption amount. As a result, if the agent consumes more at time $t$, her future selves will consume less because the remaining budget in the future becomes smaller; consequently, the budget constraint is still satisfied. To elaborate, suppose the agent's future selves' strategies are to consume $\hat k_s$ fractional of wealth at time $s$, $s=t+1,\dots, T$ with $\hat k_s\in [0,1],s=t+1,\dots, T-1$ and $\hat k_T=1$. Then, given that the agent at time $t$ consumes any amount $C_t\in [0,x_t]$, the agent's consumption in the future is $C_s = \hat k_s x_s$, $s=t+1,\dots, T$, where $x_s = x_{s-1}-C_{s-1}$, $s=t+1,\dots, T$. As a result, the aggregate consumption from time $t$ to the end is $\sum_{s=t}^TC_s = x_t$. Recall that the aggregate consumption strictly prior to time $t$ is $x_0-x_t$; so the aggregate consumption throughout the entire horizon is $x_0$ satisfying the budget constraint. Thus, at each time $t$, the agent can consume any amount up to his wealth level at that time and her future selves will adjust their consumption according to a given strategy so that the budget constraint is still satisfied.

Finally, we establish a connection between closed-loop and open-loop intra-personal equilibria. If a closed-loop equilibrium $\hat {\mathbf{u}}$ is independent of the state variable $x$, then it follows from the definition that it is also an open-loop equilibrium. For a general closed-loop equilibrium $\hat {\mathbf{u}}$, we can consider the following controlled state process:
\begin{align*}
d\hat X^{v}(s)=\hat \mu(s, \hat X^{v}(s), v(s) )ds+\hat \sigma(s,\hat X^{v}(s), v(s) )dW(s),\; s\in[t,T];\; X^{v}(t)=x,
 \end{align*}
where $\hat \mu(s,y,v):=\mu(s,y,\hat{\mathbf{u}}(s,y)+v)$, $\hat \sigma(s,y,v) : =\sigma(s,y,\hat{\mathbf{u}}(s,y)+v)$, and $v(\cdot)$ is a progressively measurable control process. We further consider the following objective function:
\begin{align*}
  \hat J(t,x;v(\cdot)):=\expect_{t,x}\left[\int_t^T \hat C\big(t,x, s, \hat X^{v}(s), v(s) )\big)ds+F\big(t, x, \hat X^{v}(T) \big)\right]\notag \\
 +G\big(t, x, \expect_{t,x} [  \hat X^{v}(T) ]  \big),
\end{align*}
where $\hat C(t,x,s,y,v):=C(t,x,s,y,\hat{\mathbf{u}}(s,y)+v)$. Then, by definition, $\hat {\mathbf{u}}$ is a closed-loop equilibrium if and only if $\hat v(\cdot)\equiv 0$ is an open-loop  equilibrium for the problem of maximizing $\hat J(t,x;v(\cdot))$ in $v(\cdot)$ with the controlled state process $\hat X^{v}$. In particular, we can characterize $\hat {\mathbf{u}}$ by a flow of forward-backward SDEs by applying the spike variation technique. In order to apply this technique, however, we need to assume that $\hat \mu(s,y,v)$ and $\hat \sigma(s,y,v)$ to be twice differentiable in $y$, which in turn requires $\hat {\mathbf{u}}$ to be twice differentiable; see \citet{YongJZhouXY:99sc} for the detailed regularity conditions needed for the spike variation technique. Thus, the spike variation technique does not seem to be advantageous over the approached reviewed in Section \ref{sec:ExHJB}.

\section{Optimal Stopping}\label{sec:Stopping}
An optimal stopping problem is one to search an optimal random time $\tau$ to stop a given, {\it uncontrollable} process $(X_t)_{t\ge 0}$ (taking values in a state space $\bbX$) in the set of stopping times with respect to the filtration generated by the process. It is well known that if the objective function of the optimal stopping problem depends on the path of $(X_t)_{t\ge 0}$ up to the stopping time only, this problem can be ``embedded" into a general control problem with (i) a closed-loop control strategy $\mathbf{u}$ taking binary values 0 and 1 representing the action of stopping and not stopping $(X_t)_{t\ge 0}$ respectively;  and (ii) a controlled state process $(\tilde X^{\mathbf{u}})_{t\ge 0}$ that is set to be $(X_t)_{t\ge 0}$ until the first time the control path under $\mathbf{u}$ takes value 0 and is set to be an {\em absorbing state} afterwards; see for instance Section 3.4 of \citet{Bertsekas2017:DynamicProgrammingOptimalControl}. We call the control strategy $\mathbf{u}$ associated with a stopping time $\tau$ in the above embedding a {\em stopping rule}, which maps each pair of time $t$ and a path of the process $X$ up to time $t$ to $\{0,1\}$. A stopping time $\tau$ is {\em Markovian} if the associated stopping rule is Markovian, i.e., it is a mapping from the time--state space  to $\{0,1\}$. With a Markovian stopping time, at each time $t$, given that the process has not yet been stopped, whether to stop at $t$ depends on the value of the process at $t$ only.


 In view of the above embedding, intra-personal equilibrium stopping rules can be defined naturally for time-inconsistent stopping problems; see for instance \citet{TWZ}, \citet{Christensen2018finding}, \citet{EbertEtal2017:Discounting}, and \citet{ChristensenLindensjo2020:OnTimeInconsistentStopping}.
In particular, \citet{TWZ} show that the smooth pasting principle, which is the main approach used to construct explicit solutions for classical time-consistent optimal stopping, may fail to find an equilibrium when one changes merely the exponential discounting to non-exponential one while keeping everything else the same. The authors also construct an explicit example in which no equilibrium exists.
These results caution blindly extending the classical approach for time-consistent
stopping to their time-inconsistent counterpart.

By exploiting special structures of stopping problems in continuous time, \citet{HuangNguyenHuu2018:TimeConsistent} propose an alternative approach to defining the optimal stopping rule for a sophisticated agen; see also applications of this approach in \citet{HuangEtal2017:StoppingBehaviors}, \citet{EbertStrack2016:NeverEverGettingStarted}, and \citet{HuangYu2021:OptimalStopping}. Precisely, consider a Markov state process
\begin{align*}
  dX_t = \mu(t,X_t)dt + \sigma(t,X_t)dW_t
\end{align*}
in $\mathbb{R}^n$, where $(W_t)_{t\ge 0}$ is an $d$-dimensional standard Brownian motion and $\mu$ and $\sigma$ are functions of time $t$ and state $x$ taking values in $\mathbb{R}^n$ and $\mathbb{R}^{n\times d}$, respectively. Following the settings in the above papers, we consider Markovian stopping times only in the following presentation, but the case of non-Markovian stopping times can be investigated  similarly.
At each time $t$ with state $x$, give that the state process has not been stopped, the agent's goal is to choose a Markovian stopping time $\tau\in [t,T]$
to maximize an objective value $J(t,x;\tau)$. Here, $J(t,x;\tau)$ can be of the form $\expect_{t,x}\left[\int_t^\tau g(t,x,s,X_s)ds + h(t,x,\tau, X_\tau)\right]$ for some functions $g$ and $h$, or be a functional of the distribution of $X_\tau$ conditional on $X_t= x$.


Recall the embedding of optimal stopping problems into a general control framework and the stopping rule associated with each stopping time as discussed at the beginning of the present subsection. With a slight abuse of notation, we use $\tau$ to denote both a stopping time and a stopping rule.
Let us now consider a given stopping rule $\tau$ and the current time--state pair $(t,x)$. If the agent decides to stop, then she has the immediate reward $J(t,x;t)$. If the agent decides not to stop at $t$ but expects her future selves will still follow the original rule $\tau$, then she will stop at time ${\cal L}^*\tau$, the first time $s>t$ at which $\tau$ would stop the process. In this case the objective  value is $J(t,x;{\cal L}^*\tau)$. Then, the optimal action of the agent at time $t$ with state $x$ is to stop if $J(t,x;t)>J(t,x;{\cal L}^*\tau)$, to continue if $J(t,x;t)<J(t,x;{\cal L}^*\tau)$, and to follow the originally assigned stopping rule $\tau$ in the break-even case $J(t,x;t)=J(t,x;{\cal L}^*\tau)$. The above plan across all time $t$ and state $x$ constitutes a {\it new} stopping rule, denoted as $\Theta \tau$, which can be proved to be feasible in the sense that it can generate stopping times; see \citet{HuangNguyenHuu2018:TimeConsistent} and \citet{HuangEtal2017:StoppingBehaviors}.

The above game-theoretic thinking shows that for any arbitrarily given stopping rule $\tau$, at any time $t$ with any state $x$, the agent finds $\Theta\tau$ to be always no worse than $\tau$, {\it assuming} that her future selves will follow $\tau$. Hence, an equilibrium stopping rule $\tau$ can be defined as one that can not be strictly improved by taking $\Theta\tau$ instead. Following \citet{BayraktarEtal2019:OnTheNotions}, we name it as a {\em mild intra-personal equilibrium} stopping rule:
\begin{definition}
\label{de:EquiStoppingIter}
  A stopping rule $\tau$ is a mild intra-personal equilibrium if $\Theta \tau=\tau$.
\end{definition}

So a mild intra-personal equilibrium is a fix-point of the operator $\Theta$.
If $\tau$ is to stop the process at any time and with any state, then it is straightforward to see that ${\cal L}^*\tau=\tau$. Consequently, by definition $\Theta \tau=\tau$ and thus $\tau$ is a mild intra-personal equilibrium. In other words, following Definition \ref{de:EquiStoppingIter}, immediate stop is {\it automatically} a (trivial) mild intra-personal equilibrium.

For a general stopping rule $\tau$, consider any time $t$ and state $x$ in the interior of the stopping region of $\tau$,  where the stoping region refers to the set of time-state pairs at which the stopping rule $\tau$ would stop the process. Then, it is also easy to see that ${\cal L}^*\tau = \tau$ at time $t$ and state $x$, so one should immediately stop under $\Theta\tau$ as well. As a result, the stopping region of $\Theta\tau$ is at least as large as that of $\tau$, if we ignore the time-state pairs that are on the boundary of the stopping region of $\tau$. Therefore, we expect the iterative sequence $\Theta^n \tau$ to converge as $n\rightarrow \infty$, and the convergent point $\tau^*$ satisfies $\tau^* = \Theta \tau^*$ and thus is a mild intra-personal equilibrium. It is, however mathematically challenging to formalize the above heuristic derivation. Rigorous proofs have been established in various settings by \citet{HuangNguyenHuu2018:TimeConsistent}, \citet{HuangEtal2017:StoppingBehaviors}, and \citet{HuangYu2021:OptimalStopping}. The above iterative algorithm, which generates a sequence $\Theta^n\tau,\; n=0,1,\dots$, not only yields a mild intra-personal equilibrium as the limit of the sequence, but also has a clear economic interpretation: each application of $\Theta$ corresponds to an additional level of strategic reasoning; see \citet{HuangNguyenHuu2018:TimeConsistent} and \citet{HuangEtal2017:StoppingBehaviors} for elaborations.

As discussed in the above, immediate stop is always a mild equilibrium; so it is expected that there exist multiple mild intra-personal equilibrium stopping rules; see \citet{HuangNguyenHuu2018:TimeConsistent} and \citet{HuangEtal2017:StoppingBehaviors}. To address the issue of multiplicity, \citet{HuangZhou2019:OptimalEquilibria} and \citet{HuangWang2020:OptimalEquilibria} consider, in the setting of an infinite-horizon, continuous-time optimal stopping under nonexponential discounting, the ``optimal" mild intra-personal equilibrium stopping rule $\tau^*$ which achieves the maximum of $J(t,x;\tau)$ over $\tau\in {\cal E}$ for all $t\in [0,T)$, $x\in \bbX$, where ${\cal E}$ is the set of all mild intra-personal equilibrium stopping rules.

 \citet{BayraktarEtal2019:OnTheNotions} compare mild intra-personal equilibrium stopping rules with weak (respectively strong) intra-personal equilibrium stopping rules obtained by embedding  optimal stopping into stochastic control  and then applying Definition \ref{de:EquilibriumFirstOrd} (respectively Definition \ref{de:EquilibriumExtd}). Assuming the objective function to be a multiplication of a discount function and a Markov process taking values in a finite or countably infinite state space, the authors prove that the optimal mild intra-personal equilibrium is a strong intra-personal equilibrium.


\section{Discretization Approach}\label{sec:Discretization}
In the discrete-time setting, an intra-personal equilibrium strategy of a sophisticated agent can be easily  defined and derived
in a backward manner starting from the last period. Thus, for a continuous-time problem, it is natural to discretize and then pass to the limit. Specifically, one partitions the continuous-time period $[0,T]$ into a finite number of subperiods, assumes the agent is able to commit in each subperiod but not beyond it, and computes the strategy chosen by the agent. Sending the length of the longest subperiod in the partition to zero, the limit of the above strategy, if it exists, can be regarded as the strategy of a sophisticated agent for the continuous-time problem. This ideas was first employed by \cite{Pollak1968:ConsistentPlanning} to study the consumption problem of \citet{Strotz1955:MyopiaInconsistency} and has recently been revisited and extensively studied by a series of papers; see for instance \cite{Yong2012:TimeInconsistent}, \cite{WeiEtal2017:TimeInconsistent}, \cite{MeiYong2019:EquilibriumStrategies}, and \cite{WangYong2019:TimeInconsistent}.

Specifically, consider the control problem in Section \ref{sec:ExHJB} and assume that in the objective function in \eqref{eq:ObjFun}, $C$ and $F$ do not depend on $x$ and $G\equiv 0$. For a partition $\Pi$ of $[0,T]$: $0=t_0<t_1<\dots <t_{N-1}<t_N=T$, we denote $\|\Pi\|:=\max_{k=1,\dots, N}|t_{k}-t_{k-1}|$. A control strategy $\hat{\mathbf{u}}^\Pi$ is an intra-personal equilibrium with respect to the partition $\Pi$ if
\begin{align}\label{eq:EquiDiscretizationCond}
  J(t_k,x_k;\hat{\mathbf{u}}^\Pi)\ge J(t_k,x_k;\mathbf{u}^\Pi_{k,\mathbf{a}})
\end{align}
for any $k=0,1,\dots, N-1$, reachable state $x_k$ at time $k$ under $\hat{\mathbf{u}}^{\Pi}$, and strategy $\mathbf{a}$, where $\mathbf{u}^\Pi_{k,\mathbf{a}}(s,\cdot) := \mathbf{a}(s,\cdot)$ for $s\in [t_k,t_{k+1})$ and $\mathbf{u}^\Pi_{k,\mathbf{a}}(s,\cdot) = \hat{\mathbf{u}}^\Pi(s,\cdot)$ for $s\in [t_{k+1},T)$. In other words, $\hat{\mathbf{u}}(s,\cdot),s\in[t_k,t_{k+1})$, is optimal for an agent who can commit in the period $[t_k,t_{k+1})$ and anticipates that her future selves will take strategy $\hat{\mathbf{u}}$ beyond time $t_{k+1}$. In the aforementioned literature, the authors define a strategy $\hat{\mathbf{u}}$ to be a {\em limiting intra-personal equilibrium} if there exists a sequence of partition $(\Pi_m)_{m\in \mathbb{N}}$ with $\lim_{m\rightarrow \infty}\|\Pi_m\|=0$ such that the state process, control process, and continuation value process under certain intra-personal equilibrium with respect to $\Pi_m$ converge to those under $\hat{\mathbf{u}}$, respectively, as $m\rightarrow \infty$. Assuming that the diffusion coefficient  of the controlled state process is independent of  control and non-degenerate and that some other conditions hold, \citet{WeiEtal2017:TimeInconsistent} prove the above convergence for any sequence of partitions with mesh size going to zero, and the limit of the continuation value function satisfies a flow of PDEs. Moreover, this flow of PDEs admits a unique solution, so the limiting intra-personal equilibrium uniquely exists. Furthermore, the limiting equilibrium is also an  equilibrium under Definition \ref{de:EquilibriumFirstOrd}.

Whether the equilibrium with respect to $\Pi$ converges when $\|\Pi\|\rightarrow 0$ for a general time-inconsistent problem, however, is still unknown. Moreover, the definition of this equilibrium relies on the assumptions that $C$ and $F$ do not depend on $x$ and  $G\equiv 0$. Otherwise, for a given partition $\Pi$, the optimal strategy  the agent at time $t_k$ implements in the subperiod $[t_k,t_{k+1})$ is {\em semi-Markovian}: the agent's action at time $s\in [t_k,t_{k+1})$ is a function of $s$, the state at $s$, and the state at $t_k$. As a result, the intra-personal equilibrium with respect to $\Pi$ is non-Markovian; so we cannot restrict  limiting equilibria to be Markov strategies.

%
%
%
%
%
%
%
%
%

\section{Applications}\label{sec:Applications}
\subsection{Present-bias Preferences}
Present-biased preferences, also known as hyperbolic discounting, refer to the following observation in intertemporal choice: when considering time preferences  between two moments, individuals become more impatient when the two moments are closer to the present time. \citet{Thaler1981:SomeEmpiricalEvidence} provides an illustrative example of present-biased preferences: some people may prefer an apple today to two apples tomorrow, but very few people would prefer an apple in a year to two apples in a year plus one day. Noted as early as in \citet{Strotz1955:MyopiaInconsistency}, present-biased preferences lead to time inconsistency. For example, consider an agent whose time preferences for having apples are as described in the above illustrative example by \citet{Thaler1981:SomeEmpiricalEvidence}. At time $0$, faced with Option A of  having one apple at time $t=365$ (days) and Option B of having two apples at time $s=366$ (days), the agent chooses Option B. When time $t=365$ arrives, however, if the agent gets to choose again, she would choose Option A. This shows that the agent in the future will change her actions planned today; hence time-inconsistency is present. For a review of the literature on present-biased preferences, see \citet{FrederickEtal2002:TimePreferences}.

In a time-separable discounted utility model, present-biased preferences can be modeled by a non-exponential discount function. For example, consider an intertemporal consumption model in continuous time for an agent. The agent's preference value of a random consumption stream $(C_s)_{s\in [t,T]}$ can be represented as
\begin{align}\label{eq:ConsumptionHBD}
  \expect_t\left[\int_t^Th(s-t)u(C_s)ds\right],
\end{align}
where $u$ is the agent's utility function, $h$ is the agent's discount function, and  $\expect_t$ denotes the expectation conditional on all the information available at time $t$. To model present-biased preferences, we assume $h(s+\Delta)/h(s)$ to be {\it strictly} increasing in $s\ge 0$ for any fixed $\Delta >0$; hence it excludes the standard exponential discount function. An example is the generalized hyperbolic discount function proposed by \citet{LoewensteinPrelec1992:Anomalies}: $h(s)= (1+\alpha s)^{-\beta/\alpha},s\ge 0$, where $\alpha>0$ and $\beta>0$ are two parameters. \citet{EbertEtal2017:Discounting} introduce a class of weighted discount functions that is broad enough to include most commonly used non-exponential discount functions in finance and economics.

In various continuous-time settings, \citet{Barro1999:Ramsey}, \citet{EkelandLazrak2006:BeingSeriousAboutNonCommitment}, \citet{EkelandILazrakA:07ti}, \citet{EkelandLazrak2010:GoldenRule}, \citet{EkelandPirvu2008:InvestConsptnWithoutCommit}, \citet{MarinSolanoNavas2010:ConsumptionPortfolio}, and \citet{EkelandEtal2012:TimeConsistent} study intra-personal equilibria for portfolio selection and consumption problems with present-biased preferences. \citet{EbertEtal2017:Discounting} and \citet{TWZ} study real option problems for agents with general weighted discount functions and derive equilibrium  investment strategies. \citet{HarrisLaibson2013:InstantaneousGratification} and \citet{GrenadierWang2007:Investment} apply a stochastic, piece-wise step discount function to a consumption problem and a real option problem, respectively, and derive intra-personal equilibrium strategies. Asset pricing for sophisticated agents with present-biased preferences and without commitment has been studied by \citet{LuttmerMariotti2003:SubjectiveDiscounting} and \citet{Bjork2017:TimeInconsistent}.

\subsection{Mean-Variance}
A popular decision criterion in finance is mean--variance, with which an agent minimizes the variance and maximizes the mean of certain random quantity, e.g., the wealth of a portfolio at the end of a period. Any mean--variance model is inherently  time inconsistent due to the variance part. To see this, consider a two-period decision problem with dates 0, 1, and 2 for an agent.  The agent is offered various options at time 1 that will yield certain payoffs at time 2. The set of options offered to the agent at time 1 depends on the outcome of a fair coin that is tossed between time 0 and 1. If the toss yields a head, the agent is offered two options at time 1: Option H1 that yields \$0 and \$200 with equal probabilities and Option H2 that yields \$50 and \$150 with equal probabilities. If the toss yields a tail, the agent is offered another two options at time 1: Option T1 that yields \$0 and \$200 with equal probabilities and Option T2 that yields \$1050 and \$1150 with equal probabilities. Suppose that at both time 0 and 1, the agent's decision criterion is to minimize the variance of the terminal payoff at time 2. At time 0, the agent has not yet observed the outcome of the toss; so she will need to make choices contingent on this outcome, i.e., she chooses between the following four plans: (H1,T1), (H1,T2), (H2,T1), and (H2,T2), where the first and second components of each of the above four plans stand for the agent's planned choice when the toss yields a head and a tail, respectively. Straightforward calculation shows that the plan (H2,T1) yields the smallest variance of the terminal payoff; so at time 0 the agent plans to choose H2 when the toss yields a head and choose T1 when the toss yields a tail. At time 1, after having observed the outcome of the toss, if the agent can choose again with the objective of minimizing the variance of the terminal payoff, she would choose H2 if the outcome is a head and T2 is the outcome is a tail. Consequently, what the agent plans at time 0 is different from what is optimal for the agent at time 1, resulting in time inconsistency.

The reason of having time inconsistency above can be seen from the following conditional variance formula: $\vari(X) = \expect[\vari(X|Y)] + \vari(\expect[X|Y])$, where $X$ stands for the terminal payoff and $Y$ denotes the outcome of the coin toss. At time 0, the agent's objective is to maximize $\vari(X)$ and at time 1, her objective is to maximize $\vari(X|Y)$. Although the plan (H2,T2) yields small variance of $X$ given the outcome of the toss $Y$ and thus a small value of the average conditional variance $\expect[\vari(X|Y)]$, it yields very different expected payoffs conditional on having a head and on having a tail, leading to a large value of $\vari(\expect[X|Y])$. Consequently, $\vari(X)$ under plan (H2,T2) is larger than under plan (H2,T1), which yields a larger value of $\expect[\vari(X|Y)]$ than the former but a much smaller value of $\vari(\expect[X|Y])$. Consequently, (H2,T1) is preferred to (H2,T2) for the agent at time 0.

A lot of recent works study intra-personal equilibrium investment strategies for agents with mean-variance preferences. For continuous-time models, see for instance \citet{BasakSChabakauri2010:DynamicMeanVariance}, \citet{BjorkEtal2011:MeanVariancewithStateDependentRiskAversion}, \citet{Pun2018:TimeConsistentMV}, \citet{BensoussanEtal2014:TimeConsistent}, \citet{CuiEtal2016:ContinuousTime}, \citet{SunEtal2016:Precommitment}, \citet{Landriault2018:EquilibriumStrategies}, \citet{BensoussanEtal2019:AParadox}, \citet{KrygerEtal2020:OptimalControl}, and \citet{HanEtal2021:RobustMV}. In all these works, the mean-variance criterion is formulated as a weighted average of the mean and variance of wealth at a terminal time, i.e., at each time $t$, the agent's objective is to maximize $\expect_t[X]-\frac{\gamma_t}{2}\vari_t(X)$, where $\expect_t$ and $\vari_t$ stand for the conditional mean and variance of the terminal wealth $X$, respectively, and $\gamma_t$ is a risk aversion parameter. Alternatively, \citet{HeJiang2017:DynamicMeanRisk} and \citet{HeJiang2020:DynamicMVFractionalKelly} study intra-personal equilibria for mean-variance investors in a constrained formulation: at each time, an investor minimizes the variance of terminal wealth with a target constraint of the expected terminal wealth. \citet{DaiEtal2017:RoboAdvising} consider a mean-variance model for log returns. \citet{HuEtal2012:TimeInconsistentStochasticLQ}, \citet{HuEtal2017:TimeInconsisten}, \citet{Czichowsky2013:TimeConsistent}, and \citet{YanWong2020:OpenLoop} investigate  open-loop intra-personal equilibria for mean-variance portfolio selection problems. For  equilibrium mean-variance insurance strategies, see for instance \citet{ZengLi2011:OptimalTimeConsistent}, \citet{LiEtal2012:OptimalTimeConsistent}, \citet{ZengEtal2013:TimeConsistentInvestment}, \citet{LiangSong2015:TimeConsistent}, and \citet{BiCai2019:OptimalInvestment}.

\subsection{Non-EUT Preferences}

There is abundant empirical and experimental evidence showing that when making choices under uncertainty, individuals do not maximize expected utility (EU); see for instance a survey by \citet{StarmerC:00neut}. Various alternatives to the EU model, which are generally referred to as {\em non-EU} models, have been proposed in the literature. Some of these models employ probability weighting functions to describe the tendency of overweighing extreme outcomes that occur with small probabilities, examples being prospect theory (PT) \citep{KahnemanDTverskyA:79pt,TverskyKahneman1992:CPT} and rank-dependent utility (RDU) theory \citep{QuigginJ:82rd}.

It has been noted that when applied to dynamic choice problems, non-EU models can lead to time inconsistency; see \citet{Machina1989:DynamicConsistency} for a review of early works discussing this issue. For illustration, consider a casino gambling problem studied by \citet{Barberis2012:Casino}: a gambler is offered 10 independent bets with equal probabilities of winning and losing \$1, plays these bets sequentially, and decides when to stop playing. Suppose at each time, the gambler's objective is to maximize the preference value of the payoff at end of the game and the preferences are represented by a non-EU model involving  a probability weighting function. We represent the cumulative payoff of playing the bets by a binomial tree with up and down movements standing for winning and losing, respectively. At time 0, the top most state (TMS) of the tree at $t=10$ represents the largest possible payoff achievable and the probability of reaching this state is extremely small ($2^{-10}$). The gambler  overweighs this state due to probability weighting and aspires  to reach it. Hence, at time 0, her plan is to play the 10-th bet if and when she has won all the previous 9 bets. Now, suppose she has played and indeed won the first 9 bets. If she has a chance to re-consider her decision of whether to play the 10-th bet at {\it that} time, she may find it no longer favorable  to play because the probability of reaching the TMS at time 10 is 1/2 and thus this state is not overweighed. Consequently, when deciding whether to play the 10-th bet conditioning on she  has won the first 9 bets, the gambler may choose differently when she is at time 0 and when she is at time 9, showing time inconsistency.

In a continuous-time, complete market, \citet{HuEtal2020:ConsistentInvestmentRDU} study a portfolio selection problem in which an agent maximizes the following RDU of her wealth $X$ at a terminal time:
\begin{align}\label{eq:RDU}
  \int_{\mathbb{R}} u(x) w(1-F_X(x)),
\end{align}
where $u$ is a utility function, $w$ is a probability weighting function, and $F_X$ is the cumulative distribution function of $X$. The authors derive an open-loop intra-personal equilibrium and show that it is in the same form as in the classical Merton model but with a properly scaled market price of risk. \citet{HeEtal2019:MedianMaximization} consider median and quantile maximization for portfolio selection, where the objective function, namely the quantile of the terminal wealth, can be regarded as a special case of RDU with a particular probability weighting function $w$. The authors study closed-loop intra-personal equilibrium and find that an affine trading strategy is an equilibrium if and only if it is a portfolio insurance strategy. \citet{EbertStrack2016:NeverEverGettingStarted} consider the optimal time to stop a diffusion process with the objective to maximize the value of the process at the stopping time under a PT model. Using the notion of mild intra-personal equilibrium as previously discussed  in Section \ref{sec:Stopping}, the authors show that under reasonable assumptions on the probability weighting functions, the only equilibrium among all two-threshold stopping rules is to immediately stop. \citet{HuangEtal2017:StoppingBehaviors} study mild intra-personal equilibrium stopping rules for an agent who wants to stop a geometric Brownian motion with the objective of maximizing the RDU value at the stopping time.

Risk measures, such as value-at-risk (VaR) and conditional value-at-risk (VaR), can also be considered to be non-EU models leading to time consistency. There are, however, few studies on intra-personal equilibria for mean-risk models in continuous time. For relevant studies in discrete-time settings, see for instance \citet{CuiEtal2019:TimeConsistent}.

Models with Knightian uncertainty or ambiguity can also result in time inconsistency. For example, the $\alpha$-maxmin model proposed by \citet{GhirardatoEtal2004:DifferentiatingAmbiguity} is dynamically inconsistent in general; see for instance \citet{BeissnerEtal2016:DynamicallyConsistent}. \citet{LiEtal2019:EquilibriumStrategies} find an open-loop intra-personal equilibrium investment strategy for an agent with $\alpha$-maximin preferences. \citet{HuangYu2021:OptimalStopping} consider a problem of stopping a one-dimensional diffusion process with preferences represented by the $\alpha$-maxmin model and study the mild intra-personal equilibrium stopping rule for the problem.


\section{Dynamically Consistent Preferences}\label{sec:DynamicConsistent}

\citet{Machina1989:DynamicConsistency} notes that, in many discussions of time inconsistency in the literature, a hidden assumption is {\em consequentialism}: at any intermediate time $t$ of a dynamic decision process, the agent employs the {\em same} preference model as used at the initial time to evaluate the choices in the {\em continuation} of the dynamic decision process from time $t$, conditional on the circumstances at time $t$. For example, consider a dynamic consumption problem for an agent with present-bias preferences and suppose that at the initial time 0, the agent's preference value for a consumption stream $(C_s)_{s\ge 0}$ is represented by $\expect[\int_0^\infty h(s)u(C_s)ds]$, where the discount function $h$ models the agent's time preferences at the initial time 0 and $u$ is the agent's utility function. The consequentialism assumption implies that at any intermediate time $t$, the agent's preferences for the continuation of the consumption stream, i.e., $(C_s)_{s\ge t}$, are represented by the same preference model as at the initial time 0, conditional on the situations at time $t$, i.e., by $\expect_t[\int_t^\infty h(s-t)u(C_s)ds]$, where the discount function $h$ and $u$ are the same as the ones in the preference model at the initial time 0. Similarly, for a dynamic choice problem with RDU preferences for the payoff at a terminal time, the consequentialism assumption stipulates that the agent uses the same utility function $u$ and probability weighting function $w$ at all intermediate times $t$ when evaluating the terminal payoff at those times. 

The consequentialism assumption, however, has not been broadly validated  because there are  few experimental or empirical studies on how individuals dynamically update their preferences.
\citet{Machina1989:DynamicConsistency} consider a class of non-EU maximizers, referred to as $\gamma$-people, who adjust their preferences dynamically over time so as to remain time consistent. The idea in \citet{Machina1989:DynamicConsistency} was further developed by \citet{KarnamEtal2016:DynamicApproaches} who propose the notion of {\em time-consistent dynamic preference models}. The idea of considering time-consistent dynamic preferences is also central in the theory of forward performance criteria proposed and developed by \citet{MusielaZariphopoulou06,MusielaZariphopoulou08,
MusielaZariphopoulou09,MusielaZariphopoulou10a,
MusielaZariphopoulou10b,MusielaZariphopoulou11}; see also \citet{HeEtal2019:ForwardRDU} for a related discussion.

Formally, consider a dynamic choice problem in a period $[0,T)$. A preference model at time 0 is specified for an agent, denoted as $J_0(u(\cdot))$, where $(u(s))_{s\in [0,T)}$ denotes the agent's dynamic choice. A family of dynamic preference models $J_t,\;t\in (0,T)$, are called time-consistent for the initial model $J_0$ if the optimal strategy under $J_0$, namely, the pre-committed strategy for the agent at time 0, is also optimal under $J_t$ for the agent at any future time $t\in (0,T)$. Note that given the pre-committed strategy at time 0, we can always find preference models at $t>0$ such that this strategy remains optimal. Thus, a more interesting question is whether we can find a family of time-consistent dynamic preference models that are of the same type as the initial preference model.

\citet{HeEtal2019:ForwardRDU} study portfolio selection in the Black-Scholes market for an agent whose initial preference model for wealth at a terminal time is represented by RDU. The authors show that there exists a family of time-consistent dynamic RDU models if and only if (i) the probability weighting function in the initial model belongs to a parametric class of functions proposed by \citet{Wang1996:PremiumCalculation}; and (ii) the parameter of the probability weighting function, the absolute risk aversion index of the utility function, and the market price of risk must be coordinated with each other over time in a specific way. \citet{CuiEtAl2012:Better}, \citet{KarnamEtal2016:DynamicApproaches}, and \citet{HeJiang2020:DynamicMVFractionalKelly} find that mean-variance models become time consistent if the dynamic trade-off between the mean and variance over time is set properly. For mean-CVaR models, where an agent maximizes the mean and minimize the CVaR at certain confidence level, \citet{PflugPichler2016:TimeInconsistent} and \citet{StrubEtal2017:DiscreteTimeMeanCVaR} note, in discrete-time settings, that time consistency is retained as long as the tradeoff between the mean and CVaR and the confidence level evolve dynamically in a certain way.

The problem of intra-personal equilibria and that of dynamically consistent preferences can be considered primal--dual to each other: 
the former finds equilibrium strategies given the time-inconsistent preferences, whereas the latter identifies preferences given the problem is time-consistent.
Diving deeper into this relationship may call for innovative mathematical analysis and result in profound economic insights.


%
%

\begin{acknowledgement}
Xue Dong He gratefully acknowledges financial support through the General Research Fund of the Research Grants Council of Hong Kong SAR (Project No. 14200917).
Xun Yu Zhou gratefully acknowledges financial support through a start-up grant at Columbia University
and through the Nie Center for Intelligent Asset Management. This author  also thanks the hospitality of  The Chinese University of Hong Kong during the summer of 2020 when the present project started.
\end{acknowledgement}
%


\bibliography{LongTitles,BibFile}

\end{document}